\immediate\write10{Package DCpic 2002/05/16 v4.0}

\catcode`!=11 

\newcount\aux%
\newcount\auxa%
\newcount\auxb%
\newcount\m%
\newcount\n%
\newcount\x%
\newcount\y%
\newcount\xl%
\newcount\yl%
\newcount\d%
\newcount\dnm%
\newcount\xa%
\newcount\xb%
\newcount\xmed%
\newcount\xc%
\newcount\xd%
\newcount\ya%
\newcount\yb%
\newcount\ymed%
\newcount\yc%
\newcount\yd
\newcount\expansao%
\newcount\tipografo
\newcount\distanciaobjmor
\newcount\tipoarco
\newif\ifpara%
\newbox\caixa%
\newbox\caixaaux%
\newif\ifnvazia%
\newif\ifvazia%
\newif\ifcompara%
\newif\ifdiferentes%
\newcount\xaux%
\newcount\yaux%
\newcount\guardaauxa%
\newcount\alt%
\newcount\larg%
\newcount\prof%
\newcount\auxqx
\newcount\auxqy
\newif\ifajusta%
\newif\ifajustadist
\def\objPartida{}%
\def\objChegada{}%
\def\objNulo{}%


\def\!vazia{:}

\def\!pilhanvazia#1{\let\arg=#1%
\if:\arg\ \nvaziafalse\vaziatrue \else \nvaziatrue\vaziafalse\fi}

\def\!coloca#1#2{\edef\pilha{#1.#2}}

\def\!guarda(#1)(#2,#3)(#4,#5,#6){\def\id{#1}%
\xaux=#2%
\yaux=#3%
\alt=#4%
\larg=#5%
\prof=#6%
}

\def\!topaux#1.#2:{\!guarda#1}
\def\!topo#1{\expandafter\!topaux#1}

\def\!popaux#1.#2:{\def\pilha{#2:}}
\def\!retira#1{\expandafter\!popaux#1}

\def\!comparaaux#1#2{\let\argA=#1\let\argB=#2%
\ifx\argA\argB\comparatrue\diferentesfalse\else\comparafalse\diferentestrue\fi}

\def\!compara#1#2{\!comparaaux{#1}{#2}}

\def\!absoluto#1#2{\n=#1%
  \ifnum \n > 0
    #2=\n
  \else
    \multiply \n by -1
    #2=\n
  \fi}

\def\solidarrow{0}

\def\solidline{2}

\def\atright{-1}
\def\atleft{1}



\def\!ajusta#1#2#3#4#5#6{\aux=#5%
  \let\auxobj=#6%
  \ifcase \tipografo    
    \ifnum\number\aux=10 
      \ajustadisttrue 
    \else
      \ajustadistfalse  
    \fi
  \else  
   \ajustadistfalse
  \fi
  \ifajustadist
   \let\pilhaaux=\pilha%
   \loop%
     \!topo{\pilha}%
     \!retira{\pilha}%
     \!compara{\id}{\auxobj}%
     \ifcompara\nvaziafalse \else\!pilhanvazia\pilha \fi%
     \ifnvazia%
   \repeat%
   \let\pilha=\pilhaaux%
   \ifvazia%
    \ifdiferentes%
     \larg=1310720
     \prof=655360%
     \alt=655360%
    \fi%
   \fi%
   \divide\larg by 131072
   \divide\prof by 65536
   \divide\alt by 65536
   \ifnum\number\y=\number\yl
    \advance\larg by 3
    \ifnum\number\larg>\aux
     #5=\larg
    \fi
   \else
    \ifnum\number\x=\number\xl
     \ifnum\number\yl>\number\y
      \ifnum\number\alt>\aux
       #5=\alt
      \fi
     \else
      \advance\prof by 5
      \ifnum\number\prof>\aux
       #5=\prof
      \fi
     \fi
    \else
     \auxqx=\x
     \advance\auxqx by -\xl
     \!absoluto{\auxqx}{\auxqx}%
     \auxqy=\y
     \advance\auxqy by -\yl
     \!absoluto{\auxqy}{\auxqy}%
     \ifnum\auxqx>\auxqy
      \ifnum\larg<10
       \larg=10
      \fi
      \advance\larg by 3
      #5=\larg
     \else
      \ifnum\yl>\y
       \ifnum\larg<10
        \larg=10
       \fi
      \advance\alt by 6
       #5=\alt
      \else
      \advance\prof by 11
       #5=\prof
      \fi
     \fi
    \fi
   \fi
\fi} 

\def\!raiz#1#2{\n=#1%
  \m=1%
  \loop
    \aux=\m%
    \advance \aux by 1%
    \multiply \aux by \aux%
    \ifnum \aux < \n%
      \advance \m by 1%
      \paratrue%
    \else\ifnum \aux=\n%
      \advance \m by 1%
      \paratrue%
       \else\parafalse%
       \fi
    \fi
  \ifpara%
  \repeat
#2=\m}

\def\!ucoord#1#2#3#4#5#6#7{\aux=#2%
  \advance \aux by -#1%
  \multiply \aux by #4%
  \divide \aux by #5%
  \ifnum #7 = -1 \multiply \aux by -1 \fi%
  \advance \aux by #3%
#6=\aux}

\def\!quadrado#1#2#3{\aux=#1%
  \advance \aux by -#2%
  \multiply \aux by \aux%
#3=\aux}

\def\!distnomemor#1#2#3#4#5#6{\setbox0=\hbox{#5}%
  \aux=#1
  \advance \aux by -#3
  \ifnum \aux=0
     \aux=\wd0 \divide \aux by 131072
     \advance \aux by 3
     #6=\aux
  \else
     \aux=#2
     \advance \aux by -#4
     \ifnum \aux=0
        \aux=\ht0 \advance \aux by \dp0 \divide \aux by 131072
        \advance \aux by 3
        #6=\aux%
     \else
     #6=3
     \fi
   \fi
}

\def\begindc#1{\!ifnextchar[{\!begindc{#1}}{\!begindc{#1}[30]}}
\def\!begindc#1[#2]{\beginpicture 
  \let\pilha=\!vazia
  \setcoordinatesystem units <1pt,1pt>
  \expansao=#2
  \ifcase #1
    \distanciaobjmor=10
    \tipoarco=0         
    \tipografo=0        
  \or
    \distanciaobjmor=2
    \tipoarco=0         
    \tipografo=1        
  \or
    \distanciaobjmor=1
    \tipoarco=2         
    \tipografo=2        
  \or
    \distanciaobjmor=8
    \tipoarco=0         
    \tipografo=3        
  \or
    \distanciaobjmor=8
    \tipoarco=2         
    \tipografo=4        
  \fi}

\def\enddc{\endpicture}

\def\mor{%
  \!ifnextchar({\!morxy}{\!morObjA}}
\def\!morxy(#1,#2){%
  \!ifnextchar({\!morxyl{#1}{#2}}{\!morObjB{#1}{#2}}}
\def\!morxyl#1#2(#3,#4){%
  \!ifnextchar[{\!mora{#1}{#2}{#3}{#4}}{\!mora{#1}{#2}{#3}{#4}[\number\distanciaobjmor,\number\distanciaobjmor]}}%
\def\!morObjA#1{%
 \let\pilhaaux=\pilha%
 \def\objPartida{#1}%
 \loop%
    \!topo\pilha%
    \!retira\pilha%
    \!compara{\id}{\objPartida}%
    \ifcompara \nvaziafalse \else \!pilhanvazia\pilha \fi%
   \ifnvazia%
 \repeat%
 \ifvazia%
  \ifdiferentes%
   Error: Incorrect label specification%
   \xaux=1%
   \yaux=1%
  \fi%
 \fi%
 \let\pilha=\pilhaaux%
 \!ifnextchar({\!morxyl{\number\xaux}{\number\yaux}}{\!morObjB{\number\xaux}{\number\yaux}}}
\def\!morObjB#1#2#3{%
  \x=#1
  \y=#2
 \def\objChegada{#3}%
 \let\pilhaaux=\pilha%
 \loop
    \!topo\pilha %
    \!retira\pilha%
    \!compara{\id}{\objChegada}%
    \ifcompara \nvaziafalse \else \!pilhanvazia\pilha \fi
   \ifnvazia
 \repeat
 \ifvazia
  \ifdiferentes%
   Error: Incorrect label specification
   \xaux=\x%
   \advance\xaux by \x%
   \yaux=\y%
   \advance\yaux by \y%
  \fi
 \fi
 \let\pilha=\pilhaaux
 \!ifnextchar[{\!mora{\number\x}{\number\y}{\number\xaux}{\number\yaux}}{\!mora{\number\x}{\number\y}{\number\xaux}{\number\yaux}[\number\distanciaobjmor,\number\distanciaobjmor]}}
\def\!mora#1#2#3#4[#5,#6]#7{%
  \!ifnextchar[{\!morb{#1}{#2}{#3}{#4}{#5}{#6}{#7}}{\!morb{#1}{#2}{#3}{#4}{#5}{#6}{#7}[1,\number\tipoarco] }}
\def\!morb#1#2#3#4#5#6#7[#8,#9]{\x=#1%
  \y=#2%
  \xl=#3%
  \yl=#4%
  \multiply \x by \expansao%
  \multiply \y by \expansao%
  \multiply \xl by \expansao%
  \multiply \yl by \expansao%
  \!quadrado{\number\x}{\number\xl}{\auxa}%
  \!quadrado{\number\y}{\number\yl}{\auxb}%
  \d=\auxa%
  \advance \d by \auxb%
  \!raiz{\d}{\d}%
  \auxa=#5
  \!compara{\objNulo}{\objPartida}%
  \ifdiferentes
   \!ajusta{\x}{\xl}{\y}{\yl}{\auxa}{\objPartida}%
   \ajustatrue
   \def\objPartida{}
  \fi
  \guardaauxa=\auxa
  \!ucoord{\number\x}{\number\xl}{\number\x}{\auxa}{\number\d}{\xa}{1}%
  \!ucoord{\number\y}{\number\yl}{\number\y}{\auxa}{\number\d}{\ya}{1}%
  \auxa=\d%
  \auxb=#6
  \!compara{\objNulo}{\objChegada}%
  \ifdiferentes
   \!ajusta{\x}{\xl}{\y}{\yl}{\auxb}{\objChegada}%
   \def\objChegada{}
  \fi
  \advance \auxa by -\auxb%
  \!ucoord{\number\x}{\number\xl}{\number\x}{\number\auxa}{\number\d}{\xb}{1}%
  \!ucoord{\number\y}{\number\yl}{\number\y}{\number\auxa}{\number\d}{\yb}{1}%
  \xmed=\xa%
  \advance \xmed by \xb%
  \divide \xmed by 2
  \ymed=\ya%
  \advance \ymed by \yb%
  \divide \ymed by 2
  \!distnomemor{\number\x}{\number\y}{\number\xl}{\number\yl}{#7}{\dnm}%
  \!ucoord{\number\y}{\number\yl}{\number\xmed}{\number\dnm}{\number\d}{\xc}{-#8}%
  \!ucoord{\number\x}{\number\xl}{\number\ymed}{\number\dnm}{\number\d}{\yc}{#8}%
\ifcase #9  
  \arrow <4pt> [.2,1.1] from {\xa} {\ya} to {\xb} {\yb}
\or  
  \setdashes
  \arrow <4pt> [.2,1.1] from {\xa} {\ya} to {\xb} {\yb}
  \setsolid
\or  
  \setlinear
  \plot {\xa} {\ya}  {\xb} {\yb} /
\or  
  \auxa=\guardaauxa
  \advance \auxa by 3%
 \!ucoord{\number\x}{\number\xl}{\number\x}{\number\auxa}{\number\d}{\xa}{1}%
 \!ucoord{\number\y}{\number\yl}{\number\y}{\number\auxa}{\number\d}{\ya}{1}%
 \!ucoord{\number\y}{\number\yl}{\number\xa}{3}{\number\d}{\xd}{-1}%
 \!ucoord{\number\x}{\number\xl}{\number\ya}{3}{\number\d}{\yd}{1}%
  \arrow <4pt> [.2,1.1] from {\xa} {\ya} to {\xb} {\yb}
  \circulararc -180 degrees from {\xa} {\ya} center at {\xd} {\yd}
\or  
  \auxa=3
 \!ucoord{\number\y}{\number\yl}{\number\xa}{\number\auxa}{\number\d}{\xmed}{-1}%
 \!ucoord{\number\x}{\number\xl}{\number\ya}{\number\auxa}{\number\d}{\ymed}{1}%
 \!ucoord{\number\y}{\number\yl}{\number\xa}{\number\auxa}{\number\d}{\xd}{1}%
 \!ucoord{\number\x}{\number\xl}{\number\ya}{\number\auxa}{\number\d}{\yd}{-1}%
  \arrow <4pt> [.2,1.1] from {\xa} {\ya} to {\xb} {\yb}
  \setlinear
  \plot {\xmed} {\ymed}  {\xd} {\yd} /
\fi
\auxa=\xl
\advance \auxa by -\x%
\ifnum \auxa=0 
  \put {#7} at {\xc} {\yc}
\else
  \auxb=\yl
  \advance \auxb by -\y%
  \ifnum \auxb=0 \put {#7} at {\xc} {\yc}
  \else 
    \ifnum \auxa > 0 
      \ifnum \auxb > 0
        \ifnum #8=1
          \put {#7} [rb] at {\xc} {\yc}
        \else 
          \put {#7} [lt] at {\xc} {\yc}
        \fi
      \else
        \ifnum #8=1
          \put {#7} [lb] at {\xc} {\yc}
        \else 
          \put {#7} [rt] at {\xc} {\yc}
        \fi
      \fi
    \else
      \ifnum \auxb > 0 
        \ifnum #8=1
          \put {#7} [rt] at {\xc} {\yc}
        \else 
          \put {#7} [lb] at {\xc} {\yc}
        \fi
      \else
        \ifnum #8=1
          \put {#7} [lt] at {\xc} {\yc}
        \else 
          \put {#7} [rb] at {\xc} {\yc}
        \fi
      \fi
    \fi
  \fi
\fi
}

\def\modifplot(#1{\!modifqcurve #1}
\def\!modifqcurve(#1,#2){\x=#1%
  \y=#2%
  \multiply \x by \expansao%
  \multiply \y by \expansao%
  \!start (\x,\y)
  \!modifQjoin}
\def\!modifQjoin(#1,#2)(#3,#4){\x=#1%
  \y=#2%
  \xl=#3%
  \yl=#4%
  \multiply \x by \expansao%
  \multiply \y by \expansao%
  \multiply \xl by \expansao%
  \multiply \yl by \expansao%
  \!qjoin (\x,\y) (\xl,\yl)             
  \!ifnextchar){\!fim}{\!modifQjoin}}
\def\!fim){\ignorespaces}

\def\setaxy(#1{\!pontosxy #1}
\def\!pontosxy(#1,#2){%
  \!maispontosxy}
\def\!maispontosxy(#1,#2)(#3,#4){%
  \!ifnextchar){\!fimxy#3,#4}{\!maispontosxy}}
\def\!fimxy#1,#2){\x=#1%
  \y=#2
  \multiply \x by \expansao
  \multiply \y by \expansao
  \xl=\x%
  \yl=\y%
  \aux=1%
  \multiply \aux by \auxa%
  \advance\xl by \aux%
  \aux=1%
  \multiply \aux by \auxb%
  \advance\yl by \aux%
  \arrow <4pt> [.2,1.1] from {\x} {\y} to {\xl} {\yl}}

\def\cmor#1 #2(#3,#4)#5{%
  \!ifnextchar[{\!cmora{#1}{#2}{#3}{#4}{#5}}{\!cmora{#1}{#2}{#3}{#4}{#5}[0] }}
\def\!cmora#1#2#3#4#5[#6]{%
  \ifcase #2
      \auxa=0
      \auxb=1
    \or
      \auxa=0
      \auxb=-1
    \or
      \auxa=1
      \auxb=0
    \or
      \auxa=-1
      \auxb=0
    \fi  
  \ifcase #6  
    \modifplot#1
    \setaxy#1
  \or  
    \setdashes
    \modifplot#1
    \setaxy#1
    \setsolid
  \or  
    \modifplot#1
  \fi  
  \x=#3%
  \y=#4%
  \multiply \x by \expansao%
  \multiply \y by \expansao%
  \put {#5} at {\x} {\y}}

\def\obj(#1,#2){%
  \!ifnextchar[{\!obja{#1}{#2}}{\!obja{#1}{#2}[Nulo]}}
\def\!obja#1#2[#3]#4{%
  \!ifnextchar[{\!objb{#1}{#2}{#3}{#4}}{\!objb{#1}{#2}{#3}{#4}[1]}}
\def\!objb#1#2#3#4[#5]{%
  \x=#1%
  \y=#2%
  \def\!pinta{\normalsize$\bullet$}
  \def\!nulo{Nulo}%
  \def\!arg{#3}%
  \!compara{\!arg}{\!nulo}%
  \ifcompara\def\!arg{#4}\fi%
  \multiply \x by \expansao%
  \multiply \y by \expansao%
  \setbox\caixa=\hbox{#4}%
  \!coloca{(\!arg)(#1,#2)(\number\ht\caixa,\number\wd\caixa,\number\dp\caixa)}{\pilha}%
  \auxa=\wd\caixa \divide \auxa by 131072 
  \advance \auxa by 5
  \auxb=\ht\caixa
  \advance \auxb by \number\dp\caixa
  \divide \auxb by 131072 
  \advance \auxb by 5
  \ifcase \tipografo    
    \put{#4} at {\x} {\y}
  \or                   
    \ifcase #5 
      \put{#4} at {\x} {\y}
    \or        
      \put{\!pinta} at {\x} {\y}
      \advance \y by \number\auxb  
      \put{#4} at {\x} {\y}
    \or        
      \put{\!pinta} at {\x} {\y}
      \advance \auxa by -2  
      \advance \auxb by -2  
      \advance \x by \number\auxa  
      \advance \y by \number\auxb  
      \put{#4} at {\x} {\y}   
    \or        
      \put{\!pinta} at {\x} {\y}
      \advance \x by \number\auxa  
      \put{#4} at {\x} {\y}   
    \or        
      \put{\!pinta} at {\x} {\y}
      \advance \auxa by -2  
      \advance \auxb by -2  
      \advance \x by \number\auxa  
      \advance \y by -\number\auxb  
      \put{#4} at {\x} {\y}   
    \or        
      \put{\!pinta} at {\x} {\y}
      \advance \y by -\number\auxb  
      \put{#4} at {\x} {\y}   
    \or        
      \put{\!pinta} at {\x} {\y}
      \advance \auxa by -2  
      \advance \auxb by -2  
      \advance \x by -\number\auxa  
      \advance \y by -\number\auxb  
      \put{#4} at {\x} {\y}   
    \or        
      \put{\!pinta} at {\x} {\y}
      \advance \x by -\number\auxa  
      \put{#4} at {\x} {\y}   
    \or        
      \put{\!pinta} at {\x} {\y}
      \advance \auxa by -2  
      \advance \auxb by -2  
      \advance \x by -\number\auxa  
      \advance \y by \number\auxb  
      \put{#4} at {\x} {\y}   
    \fi
  \or                   
    \ifcase #5 
      \put{#4} at {\x} {\y}
    \or        
      \put{\!pinta} at {\x} {\y}
      \advance \y by \number\auxb  
      \put{#4} at {\x} {\y}
    \or        
      \put{\!pinta} at {\x} {\y}
      \advance \auxa by -2  
      \advance \auxb by -2  
      \advance \x by \number\auxa  
      \advance \y by \number\auxb  
      \put{#4} at {\x} {\y}   
    \or        
      \put{\!pinta} at {\x} {\y}
      \advance \x by \number\auxa  
      \put{#4} at {\x} {\y}   
    \or        
      \put{\!pinta} at {\x} {\y}
      \advance \auxa by -2  
      \advance \auxb by -2
      \advance \x by \number\auxa  
      \advance \y by -\number\auxb 
      \put{#4} at {\x} {\y}   
    \or        
      \put{\!pinta} at {\x} {\y}
      \advance \y by -\number\auxb 
      \put{#4} at {\x} {\y}   
    \or        
      \put{\!pinta} at {\x} {\y}
      \advance \auxa by -2  
      \advance \auxb by -2
      \advance \x by -\number\auxa 
      \advance \y by -\number\auxb 
      \put{#4} at {\x} {\y}   
    \or        
      \put{\!pinta} at {\x} {\y}
      \advance \x by -\number\auxa 
      \put{#4} at {\x} {\y}   
    \or        
      \put{\!pinta} at {\x} {\y}
      \advance \auxa by -2  
      \advance \auxb by -2
      \advance \x by -\number\auxa 
      \advance \y by \number\auxb  
      \put{#4} at {\x} {\y}   
    \fi
   \else 
     \ifnum\auxa<\auxb 
       \aux=\auxb
     \else
       \aux=\auxa
     \fi
     \ifdim\wd\caixa<1em
       \dimen99 = 1em
       \aux=\dimen99 \divide \aux by 131072 
       \advance \aux by 5
     \fi
     \advance\aux by -2 
     \multiply\aux by 2 %
     \ifnum\aux<30
       \put{\circle{\aux}} [Bl] at {\x} {\y}
     \else
       \multiply\auxa by 2
       \multiply\auxb by 2
       \put{\oval(\auxa,\auxb)} [Bl] at {\x} {\y}
     \fi
     \put{#4} at {\x} {\y}
   \fi   
}

\catcode`!=12 

\tolerance=10000
\hsize=15truecm \hoffset 0.4truecm           
\vsize=21truecm \voffset 1.5truecm


\def\Abstract#1\par{\centerline{\vbox{\hsize 11truecm
    \noindent {\bf Abstract.} #1\hfil}} \vskip 0.75cm}

\newcount\numsezione
\newcount\numcapitolo
\def\Capitolo #1{\advance \numcapitolo by 1 \numsezione=0
       \penalty -1000
        \bigskip
       {\bf \the \numcapitolo. \ #1 \hfill}
       \bigskip 
       \nobreak }

\def\Sezione #1{\ifnum\numsezione>0\goodbreak\fi
\noindent
       \global\advance \numsezione by 1
       \bigskip 
       \mark{ {\bf \the\numcapitolo.\the\numsezione}\enspace #1}
       {\bf \the \numcapitolo.\the \numsezione \ #1 \hfill}
       \nobreak\bigskip}

\def\References{\vskip1cm\goodbreak\noindent
    {\bf References \hfill}
    \vskip 0.25cm \nobreak}

\newcount \numproclaim

 \long\def\proclaim#1.#2 \par{ \advance\numproclaim by 1
         \goodbreak\ifdim\lastskip<\bigskipamount\removelastskip\bigskip\fi
         \noindent
        {\bf #1 \the \numproclaim.}
        {\sl #2\par}
         \bigskip}

\newdimen\myitemindent    

\myitemindent = 20 pt     
\def\myitem#1{\par\leavevmode\hbox to \myitemindent
               {\hbox{#1} \hfil}\ignorespaces}

\def\bull{{\vrule height.9ex width.8ex depth-.1ex}}

\def\k#1{{\cal #1}}

\def\lra{\longrightarrow}

\def\u#1{{\underline#1}}

\def\PP{\vbox{\hbox to 8.9pt{I\hskip-2.1pt P\hfill}}}
\def\II{\vbox{\hbox to 8.9pt{I\hskip-2.1pt I\hfill}}}

\font\tenmsbm=msbm10
\font\sevenmsbm=msbm7
\font\fivemsbm=msbm5
\newfam\amsBfam
\textfont\amsBfam=\tenmsbm
\scriptfont\amsBfam=\sevenmsbm
\scriptscriptfont\amsBfam=\fivemsbm
\def\bbb{\fam\the\amsBfam\tenmsbm}
\def\Bbb#1{{\bf #1}}

\font\tenmsam=msam10
\font\sevenmsam=msam7
\font\fivemsam=msam5
\newfam\amsAfam
\textfont\amsAfam=\tenmsam
\scriptfont\amsAfam=\sevenmsam
\scriptscriptfont\amsAfam=\fivemsam
\def\aaa{\fam\the\amsAfam\tenmsam}

\mathchardef\gul"3\the\amsAfam52
\mathchardef\lug"3\the\amsAfam51
\mathchardef\R"5\the\amsBfam52

		\catcode`\"=12
	\font\kropa=lcircle10 scaled 1700
	\def\ybl{\setbox0=\hbox{\kropa \char"70} \kern1.5pt \raise.35pt \box0}
		\catcode`\"=\active

\input pictex

\def\bullet{}
	\def\normalsize{}

\catcode`\"=12
	\font\syy=cmsy10 scaled 1150
	\font\syx=cmsy10 scaled 1650

	\def\PROD{{\setbox0=\hbox{\syy\char"02}\box0}}
	\def\POUT{{\setbox0=\hbox{\syx\char"0E}\box0}}

\catcode`\"=\active

	\def\lra{\longrightarrow }

	\def\u#1{{\underline{#1}}}

	\def\prod#1{\mathop{\POUT}\limits_{{\,#1\,}}}
	
	\def\pback#1{\mathop{\PROD}\limits_{{\,#1\,}}}
	
	\def\pout#1{\mathop{\PROD}\limits^{{\,#1\,}}}

\numproclaim=0
\parindent=0pt

\centerline {\bf{The category of local algebras and points proches.}}
\vskip 2\baselineskip 
\centerline {MARGHERITA BARILE$^1$}
\centerline{FIORELLA BARONE$^2$}
\centerline{WLODZIMIERZ M. TULCZYJEW$^3$\footnote{(*)}{Supported by PRIN SINTESI.}}
\centerline{\hphantom{X}}
\centerline {$^1$Dipartimento di Matematica, Via Orabona, 4}
\centerline {Universit\`a di Bari, 70125 Bari, Italy}
\centerline {barile@dm.uniba.it}
\centerline {Tel. 0039 0805442711 -- Fax 0039 0805963612} 
\centerline{\hphantom{X}}
\centerline {$^2$Dipartimento di Matematica, Via Orabona, 4}
\centerline {Universit\`a di Bari, 70125 Bari, Italy}
\centerline {barone@dm.uniba.it} 
\centerline {Tel. 0039 0805442711 -- Fax 0039 0805963612}
\centerline{\hphantom{X}}
\centerline {$^3$Dipartimento di Fisica} 
\centerline {Universit\`a di Camerino,  Camerino, Italy}
\centerline {tulczy@libero.it} 
\centerline {Tel. 0039 0737519748 -- Fax 0039 0737402529}

\vskip 3\baselineskip
\Abstract 
Categorial methods for generating new local algebras from old ones are presented.
A direct proof of the differential structure of the prolongations of a manifold is proposed.

\vskip 2\baselineskip

\vskip 1\baselineskip
\noindent A.M.S. CLASSIFICATION: 13H10, 58A20

\Capitolo {Introduction.}

In his paper on the theory of points proches [15], A. Weil, generalizing the notion of jet prolongation due to C.  Ehresmann [3], introduced the idea that to any local algebra $A$ corresponds a covariant functor which associates  with a differential manifold $P$ a fibre bundle $AP$ over $P$. The bundle $AP$ is called the prolongation of $P$ of kind $A$ and its elements the $A$-points, or points proches, on $P$.

From this paper, several lines of study followed. Among them, the study of the structures one can lift on the prolongations of a manifold ([4], [11], [12]); the theory of product preserving, or natural bundles ([2], [5], [6], [7], [8]); the study of differential invariants on sheaves of tangent vector fields ([13], [14]).

Since the starting point is, in any case, a local algebra, it seemed of some interest to study the category of local algebras and, in particular, the categorial methods for generating new local algebras from old ones.
This is the content of the first section of the paper.

It is known that each local algebra is isomorphic to an algebra of (generalized) truncated polynomials. Thus, in the second section of the paper, we study these kind of local algebras, performing the constructions introduced in the previous part and finding, in particular, examples of local algebras which are isomorphic as $A$-modules but not as algebras.

Finally, in the third section, we study the differential structure of the prolongation bundles.
Weil stated the existence of this structure, and it has been explicitly shown in [7]. In this paper we propose an alternative proof, which is more direct and concise. 

We consider paracompact manifolds modelled over finite dimensional affine spaces.
We make use of the notion of polynomial function on an affine space and of some related results presented in [1], which we recall, without proofs, in the Appendix. 

We refer to [10] for all basic notions of category theory we use.

\Capitolo {The Category of Local Algebras.}

We will study the category of local algebras and, in particular, the categorial methods for generating new local algebras from old ones.
For the sake of completeness, we will first recall the main definitions and basic properties.

\Sezione {Objects and Morphisms.}

Let ${\cal C}$ be the category of associative, commutative, real  algebras  and algebra homomorphisms.

We will refer to them as {\it algebras} and {\it  morphisms} respectively.

\proclaim Definition.
A {\it local algebra} $A$ is an algebra with identity such that 
\myitem{(i)} $A$ is a finite dimensional vector space over ${\R}$;
\myitem{(ii)} $A$ has a  unique maximal ideal $I_A$ and $A/I_A\simeq{\R}$.
\hfill$\bull$

A morphism which preserves the identity will be called a  {\it local algebra morphism}.

We denote by ${\cal A}$ the subcategory of ${\cal C}$ of  local algebras and local algebra morphisms. 

Unless otherwise noted, all commutative diagrams considered in this section will be assumed to be local algebra diagrams.

\vskip 0.5\baselineskip
As a consequence of condition (ii) of Definition 1, we have an epimorphism of local algebras
$$0_A : A \to \R$$
whose kernel is $I_A$.
If $1_A$ is the identity of $A$, it follows that $A = {\R}1_A + I_A$. 
The subalgebra $\R 1_A$ of $A$ is mapped isomorphically onto $\R$ by $0_A$, via  $\lambda 1_A\mapsto \lambda$.
This gives a natural way to identify $\R 1_A$ with $\R$.  
Then we have that 
$$A = \R + I_A\,,\eqno(1)$$
where the sum is direct.
This means that for every element $a$ of $A$ there is a unique decomposition $a=\alpha+\u{a}$, where $\alpha\in\R$ and $\u{a}\in I_A$; $\alpha$ will be called the {\it finite part} of $a$.
Conversely, for every object $J$ of ${\cal C}$ without identity, we can endow the direct product 
$$A=\R\times J \eqno(2)$$
of real vector spaces with a local algebra structure by defining the product: $(\alpha_1,\u a_1)(\alpha_2,\u a_2)=(\alpha_1\alpha_2, \alpha_1\u a_2 +\alpha_2 \u a_1)$.

\vskip 0.5\baselineskip
Condition (i) of Definition 1 implies that all ideals of $A$ are finitely generated as vector spaces over $\R$.
A fortiori they are finitely generated as modules over $A$.
Hence  $A$ is a Noetherian ring.
Furthermore,  finite dimension over $\R$ forces the descending sequence of powers of $I_A$ to be stationary.
By Nakayama's Lemma (cf. [9], Th.2.2),  this implies that  $I_A^{\ell+1}=(0)$ for some integer $\ell$.
The least integer $\ell$ with this property is called the {\it height} of $A$.

\vskip0.5\baselineskip
Every local algebra morphism $\varphi: A\to B$ maps $\R$ identically onto $\R$.
In particular, $\varphi$ preserves the finite part of every element; therefore, for all $\u{a}\in I_A$, we have that $\varphi(\u{a})\in I_B$, i.e., $\varphi(I_A)\subset I_B$.
It follows that $0_A$ is the only local algebra morphism from $A$ to $\R$, and that the inclusion mapping of $\R$ into $A$ is the only local algebra morphism from $\R$ to $A$. 
Hence, in the category  ${\cal A}$, $\R$ is both the final and the initial object, i.e., it is the zero object of ${\cal A}$.

\vskip0.5\baselineskip
\Sezione {Subobjects and Quotient Objects.}

(i)\enspace A local algebra $B$ is a subobject of the local algebra $A$ in ${\cal A}$ if there is a local algebra monomorphism $\iota:B\to A$.
Up to an isomorphism, $B$ can be identified with $\iota(B)=\R +\iota(I_B)$.
Since $\iota(I_B)$ is a subspace of $I_A$ which is closed under multiplication,  $\iota(B)$ is an object
of ${\cal A}$. We will call it a {\it local subalgebra} of $A$.
Conversely, every subspace $J$ of $I_A$ which is closed under multiplication gives rise to a local subalgebra
$\R+J$ of $A$, and this correspondence defines  all subalgebras of $A$.

The above arguments show that for every local algebra morphism $\varphi:A\to A'$,
$\varphi(A)=\R+\varphi(I_A)$ is a local subalgebra of $A'$: it is the image of $\varphi$ in both categories
${\cal C}$ and ${\cal A}$.

Similarly one sees that $\R+\varphi^{-1}(I_{A'})$ is the inverse image of $\varphi$ in both categories.

On the contrary, the kernel of $\varphi$ in ${\cal C}$ ($\ker\varphi$) and the kernel of $\varphi$ in ${\cal A}$ (Ker$\,\varphi$) are distinct, since 
$${\rm Ker}\,\varphi=\R+\ker\varphi\;.$$
We will call {\it normal subalgebra} of $A$ any subalgebra $\R+J$, where $J$ is a proper ideal of $A$. 
The normal subalgebras of $A$ are the kernels in ${\cal A}$ of all local algebra morphisms defined on $A$.

\vskip0.5\baselineskip
(ii)\enspace A quotient in ${\cal A}$ of the local algebra $A$ is a local algebra $C$ for which there is a
local algebra epimorphism $\pi:A\to C$. 
Then $C$ is isomorphic to the quotient algebra $A/\ker\pi=\R+I_A/\ker\pi$, which is a local algebra with maximal ideal $I_A/\ker\pi$.

This allows us to  identify the quotient objects of $A$ in ${\cal A}$ with the local algebras of the type  
$A/J=\R+I_A/J$, where $J$ is any proper ideal  of $A$.

\vskip0.5\baselineskip
Given a local algebra morphism $\varphi: A\lra A'$, the cokernel  of $\varphi$ in both 
${\cal C}$ and ${\cal A}$ is
$${\rm Cok}\,\varphi = {\R}+I_{A'}/\bigcap_{J\ {\rm proper\ ideal},\atop \varphi(I_A)\subset J}J\;.$$

\Sezione {Products and Coproducts.}

Let $A_1$ and $A_2$ be local algebras.

\vskip0.5\baselineskip
(i) \enspace The usual product in $\k C$ does not yield a product in $\k A$, because the direct product $A_1\times A_2$ is not a local algebra.
In fact it has two maximal ideals, $A_1\times I_{A_2}$ and $I_{A_1}\times A_2$.

If we consider the subalgebra without identity  $I_{A_1}\times I_{A_2}$, following (2) we can construct the local algebra 
$$\R\times I_{A_1}\times I_{A_2}$$
which is isomorphic to the following subalgebra of $A_1\times A_2$:
$$\R(1_{A_1},1_{A_2}) + I_{A_1}\times I_{A_2}\,.$$
We denote it by
$$A_1\times_{\R}A_2:=\R + I_{A_1}\times I_{A_2}\,,$$
where we have used  the convention introduced in (1).

We can now define a functor $\times_{\R}$ which associates with every pair of local algebras $A_1$ and $A_2$ the  diagram
$$\vcenter{
\begindc{0}[13]
\obj(1,4){$A_1$}
\obj(9,6){$A_2$}
\obj(6,9){$A_1\times_{\R}A_2$}
\mor(6,9)(9,6){$Pr_{A_2}$}
\mor(6,9)(1,4){$Pr_{A_1}$}[\atright,\solidarrow]
\enddc}
\eqno(3)$$
where 
$$\eqalign{
Pr_{A_i}: A_1\times_{\R}A_2&\lra A_i\cr
\alpha+(\u a_1,\u a_2)&\mapsto \alpha+\u a_i\qquad\qquad (i=1,2),
}$$
are the natural projections, and with every pair of local algebra morphisms $\xi_1 : A_1 \lra B_1$ and $ \xi_1 : A_2 \lra B_2$ the diagram

\vskip0.5\baselineskip
\centerline{
\begindc{0}[13]
\obj(1,4){$A_1$}
\obj(14,4){$B_1$}
\obj(9,6){$A_2$}
\obj(22,6){$B_2$}
\obj(6,9){$A_1\times_{\R}A_2$}
\obj(19,9){$B_1\times_{\R}B_2$}
\mor(1,4)(14,4){$\phantom{bbbb}\xi_1$}
\mor(9,6)(22,6){$\xi_2\phantom{bbbbbbbbbbbb}$}
\mor(7,9)(18,9){$\xi_1\times_{\R}\xi_2$}
\mor(6,9)(9,6){$Pr_{A_2}$}
\mor(19,9)(22,6){$Pr_{B_2}$}
\mor(6,9)(1,4){$Pr_{A_1}$}[\atright,\solidarrow]
\mor(19,9)(14,4){$Pr_{B_1}$}[\atright,\solidarrow]
\enddc}

where, for any $\alpha+(\u{a}_1,\u{a}_2) \in A_1\times_{\R}A_2$,
$$\xi_1\times_{\R}\xi_2 \big(\alpha+(\u{a}_1,\u{a}_2)\big) = \alpha+\big(\xi_1(\u{a}_1),\xi_2(\u{a}_2)\big).$$
Since diagram (3) fulfills the universal property of products, we will call $\times_{\R}$ the {\it product functor} in 
$\k A$ and we will refer to $A_1\times_{\R}A_2$ as the product of $A_1$ and $A_2$.

Generalizing the above construction, we can consider the product of any finite number of local algebras.

The associative property is trivially fulfilled.

\vskip0.5\baselineskip

(ii) \enspace Similarly, the {\it coproduct functor} associates with $A_1$, $A_2$ the  diagram:
$$\vcenter{
\begindc{0}[13]
\obj(1,4){$A_1$}
\obj(9,6){$A_2$}
\obj(6,9){$A_1\times_{\R}A_2$}
\mor(9,6)(6,9){$In_{A_2}$}[\atright,\solidarrow]
\mor(1,4)(6,9){$In_{A_1}$}
\enddc}
\eqno(4)$$
where
$$\eqalign{
In_{A_1}: A_1&\lra A_1\times_{\R}A_2\cr
\alpha+\u a_1 &\mapsto \alpha+(\u a_1,0)\;,\cr
}$$
and
$$\eqalign{
In_{A_2}: A_2&\lra A_1\times_{\R}A_2\cr
\alpha+\u a_2 &\mapsto \alpha+(0,\u a_2)\cr
}$$
are the natural injections.
It acts in an obvious way on pairs of local algebra morphisms.

\vskip0.5\baselineskip
It is easy to verify that diagrams (3) and (4) determine a biproduct of $A_1$ and $A_2$ in ${\cal A}$.

\Sezione {Relative Products.}

\proclaim Definition.
Let $\pi_{A_1}:A_1\to B$ and $\pi_{A_2}:A_2\to B$ be local algebra epimorphisms. 
The following triple product of local algebras
$$\eqalign{
A_1 \prod{(\pi_{A_1},\pi_{A_2})}A_2 
&:= B\times_{\R}{\rm Ker}\,\pi_{A_1}\times_{\R}{\rm Ker}\,\pi_{A_2}\cr
&\phantom{ : }= \R + I_B\times \ker \pi_{A_1} \times \ker \pi_{A_2} \cr
}$$
will be called the {\it product} of $A_1$  and $A_2$ relative to  $\pi_{A_1}$ and $\pi_{A_2}$.
\hfill$\bull$

Note that, since $\ker 0_{A_i} = I_{A_i}$ for $i=1,2$, we have that
$$A_1 \prod{(0_{A_1},0_{A_2})}A_2 = A_1 \times_{\R} A_2\, .$$
This shows that Definition 2 generalizes the product introduced in 2.3.

\vskip0.5\baselineskip
We can now define a functor $\prod{(\cdot,\cdot)}$ which associates with every diagram of epimorphisms

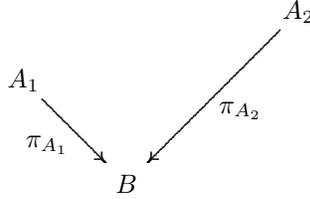
\centerline{
\begindc{0}[13]
\obj(4,1){$B$}
\obj(1,4){$A_1$}
\obj(9,6){$A_2$}
\mor(1,4)(4,1){$\pi_{A_1}$}[\atright,\solidarrow]
\mor(9,6)(4,1){$\pi_{A_2}$}
\enddc}
the diagram

$$\vcenter{
\begindc{0}[13]
\obj(4,1){$B$}
\obj(4,6){$A_1 \prod{(\pi_{A_1},\pi_{A_2})}A_2$}
\mor(4,6)(4,1){$Pr_B$}[\atright,\solidarrow]
\enddc}
\eqno(5)$$

and with every diagram

\centerline{
\begindc{0}[13]
\obj(4,1){$B$}
\obj(17,1){$D$}
\obj(1,4){$A_1$}
\obj(14,4){$C_1$}
\obj(9,6){$A_2$}
\obj(22,6){$C_2$}
\mor(4,1)(17,1){$\eta$}
\mor(1,4)(14,4){$\phantom{bbbbbbbbb}\xi_1$}
\mor(9,6)(22,6){$\xi_2$}
\mor(1,4)(4,1){$\pi_{A_1}$}[\atright,\solidarrow]
\mor(14,4)(17,1){$\pi_{C_1}$}[\atright,\solidarrow]
\mor(9,6)(4,1){$\pi_{A_2}$}
\mor(22,6)(17,1){$\pi_{C_2}$}
\enddc}

the diagram

\centerline{
\begindc{0}[13]
\obj(1,1){$B$}
\obj(10,1){$D$}
\obj(1,6){$A_1 \prod{(\pi_{A_1},\pi_{A_2})}A_2$}
\obj(10,6){$C_1 \prod{(\pi_{C_1},\pi_{C_2})}C_2$}
\mor(3,6)(8,6){${\xi_1\prod{\eta}\xi_2}$}[\atleft,\solidarrow ]
\mor(1,1)(10,1){$\eta$}[\atright,\solidarrow]
\mor(1,6)(1,1){$Pr_B$}[\atright,\solidarrow]
\mor(10,6)(10,1){$Pr_D$}
\enddc}

where, for any $\alpha+(\u{b},\u{a}_1,\u{a}_2) \in A_1 \prod{(\pi_{A_1},\pi_{A_2})}A_2$,
$$\xi_1\prod{}\xi_2 \big(\alpha+(\u{b},\u{a}_1,\u{a}_2)\big) =
\alpha+\big(\eta(\u{b}),\xi_1(\u{a}_1),\xi_2(\u{a}_2)\big).$$

Definition 2 can be extended to any finite number of local algebra epimorphisms.
For example, for $n=3$, consider the product
$$B\times_{\R}{\rm Ker}\,\pi_{A_1}\times_{\R}{\rm Ker}\,\pi_{A_2}\times_{\R}{\rm Ker}\,\pi_{A_3}\,.\eqno(6)$$
In view of diagram (5), we can also consider the product
$$\left(A_1 \prod{(\pi_{A_1},\pi_{A_2})}A_2\right)  \prod{(Pr_B,\pi_{A_3})} A_3 \,.\eqno(7)$$
Since
$${\rm Ker}\,Pr_B = {\rm Ker}\,\pi_{A_1}\times_{\R}{\rm Ker}\,\pi_{A_2}$$
and $\times_{\R}$ is associative, we have that (6) and (7) are equal. 
Similarly one can prove that (6) is also equal to
$$A_1 \prod{(\pi_{A_1},Pr_B)} \left(A_2  \prod{(\pi_{A_2},\pi_{A_3})} A_3\right)\,.$$
This shows the associativity of functor $ \prod{(\cdot,\cdot)}$.

\Sezione {Pullbacks and Pushouts.}

(i)\enspace Let $\varphi_{A_1}:A_1\to B$ and $\varphi_{A_2}:A_2\to B$ be local algebra morphisms.

Let us introduce the local algebra
$$A_1 \pback{(\varphi_{A_1},\varphi_{A_2})}A_2 =
{\R}+\{(\u a_1,\u a_2)\in I_{A_1}\times I_{A_2}\mid \varphi_{A_1}(\u a_1) = \varphi_{A_2}(\u a_2)\} .\eqno(8)$$
and notice that it is isomorphic to the following subalgebra of  $A_1\times A_2\,$:
$$\{ (a_1,a_2)\mid \varphi_{A_1}(a_1) = \varphi_{A_2}(a_2) \} = \{ (\alpha_1 + \u a_1,\alpha_2 + \u a_2) \mid \alpha_1
= \alpha_2 \ {\rm and } \ \varphi_{A_1}(\u a_1) = \varphi_{A_2}(\u a_2) \}.$$

We can now define a functor $\pback{(\cdot,\cdot)}$ which associates with every diagram

$$\vcenter
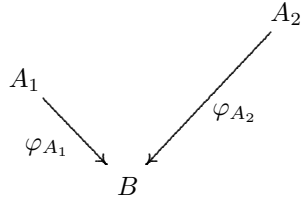
{
\begindc{0}[1]
\obj(40,10){$B$}
\obj(1,50){$A_1$}
\obj(100,75){$A_2$}
\mor(1,50)(40,10){$\varphi_{A_1}$}[\atright,\solidarrow]
\mor(100,75)(40,10){$\varphi_{A_2}$}
\enddc}
$$

the diagram

$$\vcenter{
\begindc{0}[1]
\obj(40,10){$B$}
\obj(1,50){$A_1$}
\obj(100,75){$A_2$}
\obj(60,120){$A_1\pback{(\varphi_{A_1},\varphi_{A_2})}A_2$}
\mor(1,50)(40,10){$\varphi_{A_1}$}[\atright,\solidarrow]
\mor(60,112)(100,75){$Pr_{A_2}$}
\mor(55,112)(1,50){$Pr_{A_1}$}[\atright,\solidarrow]
\mor(100,75)(40,10){$\varphi_{A_2}$}
\enddc}
\eqno(9)$$

where $Pr_{A_1}$ and $Pr_{A_2}$, are the restrictions of the natural projections defined in 2.3, and associates with every diagram

$$\vcenter
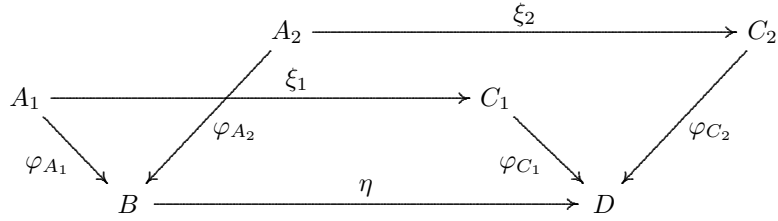
{
\begindc{0}[1]
\obj(40,10){$B$}
\obj(1,50){$A_1$}
\obj(100,75){$A_2$}
\obj(220,10){$D$}
\obj(179,50){$C_1$}
\obj(280,75){$C_2$}
\mor(1,50)(40,10){$\varphi_{A_1}$}[\atright,\solidarrow]
\mor(100,75)(40,10){$\varphi_{A_2}$}
\mor(179,50)(220,10){$\varphi_{C_1}$}[\atright,\solidarrow]
\mor(280,75)(220,10){$\varphi_{C_2}$}
\mor(1,50)(179,50){\phantom{ddddd}$\xi_1$}
\mor(100,75)(280,75){$\xi_2$}
\mor(40,10)(220,10){$\eta$}
\enddc}
$$

the diagram

$$\vcenter
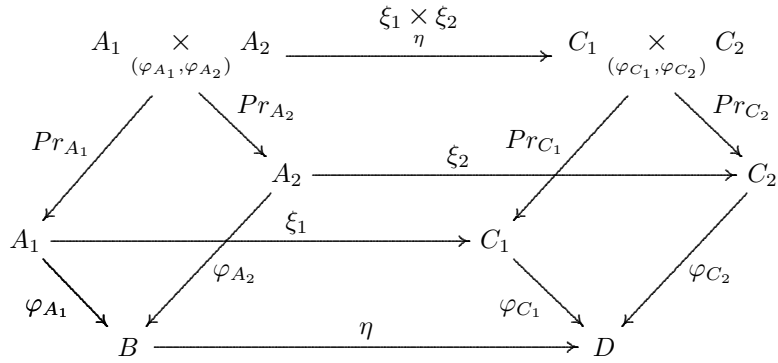
{
\begindc{0}[1]
\obj(40,10){$B$}
\obj(1,50){$A_1$}
\obj(100,75){$A_2$}
\obj(60,120){$A_1\pback{(\varphi_{A_1},\varphi_{A_2})}A_2$}
\mor(1,50)(40,10){$\varphi_{A_1}$}[\atright,\solidarrow]
\mor(60,112)(100,75){$Pr_{A_2}$}
\mor(55,112)(1,50){$Pr_{A_1}$}[\atright,\solidarrow]
\mor(100,75)(40,10){$\varphi_{A_2}$}
\obj(220,10){$D$}
\obj(179,50){$C_1$}
\obj(280,75){$C_2$}
\obj(240,120){$C_1\pback{(\varphi_{C_1},\varphi_{C_2})}C_2$}
\mor(1,50)(40,10){$\varphi_{A_1}$}[\atright,\solidarrow]
\mor(240,112)(280,75){$Pr_{C_2}$}
\mor(235,112)(179,50){$Pr_{C_1}$}[\atright,\solidarrow]
\mor(179,50)(220,10){$\varphi_{C_1}$}[\atright,\solidarrow]
\mor(280,75)(220,10){$\varphi_{C_2}$}
\mor(40,10)(220,10){$\eta$}
\mor(1,50)(179,50){\phantom{ddddd}$\xi_1$}
\mor(100,75)(280,75){$\xi_2$\phantom{ddddddddd}}
\mor(90,120)(210,120){$\xi_1\pback{\eta}\xi_2$}
\enddc}
$$

We will call $\pback{(\cdot,\cdot)}$ the {\it pullback functor} in $\k A$, since  diagram (9) fulfills the corresponding universal property (in both categories $\k A$ and $\k C$) and we will refer to $A_1\pback{(\varphi_{A_1},\varphi_{A_2})}A_2$ as the pullback of $\varphi_{A_1}$ and $\varphi_{A_2}$. 

It is possible to consider the pullback of an arbitrary number of local algebra morphisms following the usual categorial method.
The associative property is trivially fulfilled.

It is worth pointing out that $A_1\times_{{\R}}A_2 = A_1 \pback{(0_{A_1},0_{A_2})}A_2$.
This shows that the pullback gives one more generalization of the product defined in 2.3.

\vskip0.5\baselineskip

(ii)\enspace Let $\pi_{A_1}:A\to A_1$ and $\pi_{A_2}:A\to A_2$ be local algebra epimorphisms. 

Consider the quotient 
$$A/(\ker\pi_{A_1}+\ker\pi_{A_2})$$
together with the natural epimorphism
$$\pi:A\lra A/(\ker\pi_{A_1}+\ker\pi_{A_2})\,.$$
It is well known that there are unique local algebra epimorphisms 
$$Ep_{A_i} : A_i \lra A/(\ker\pi_{A_1}+\ker\pi_{A_2})\quad\quad (i=1,2)$$ 
such that
$$Ep_{A_1}\circ\pi_{A_1} = \pi = Ep_{A_2}\circ\pi_{A_2}\,.$$
We introduce the local algebra
$$A_1\pout{(\pi_{A_1},\pi_{A_2})}A_2 := A/(\ker\pi_{A_1}+\ker\pi_{A_2})\eqno(10)$$
and define a functor $\pout{(\cdot,\cdot)}$ which associates with every diagram of epimorphisms

$$\vcenter
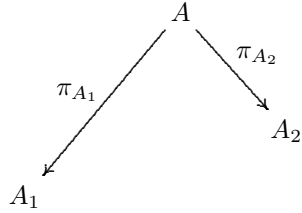
{
\begindc{0}[1]
\obj(1,50){$A_1$}
\obj(100,75){$A_2$}
\obj(60,120){$A$}
\mor(60,120)(100,75){$\pi_{A_2}$}
\mor(60,120)(1,50){$\pi_{A_1}$}[\atright,\solidarrow]
\enddc}
$$

the diagram

$$\vcenter
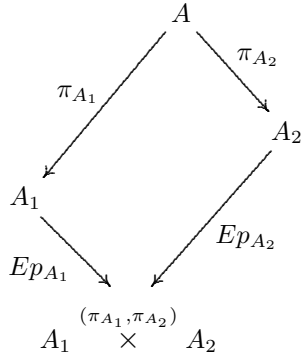
{
\begindc{0}[1]
\obj(40,1){$A_1\pout{(\pi_{A_1},\pi_{A_2})}A_2$}
\obj(1,50){$A_1$}
\obj(100,75){$A_2$}
\obj(60,120){$A$}
\mor(1,50)(40,9){$Ep_{A_1}$}[\atright,\solidarrow]
\mor(60,120)(100,75){$\pi_{A_2}$}
\mor(60,120)(1,50){$\pi_{A_1}$}[\atright,\solidarrow]
\mor(100,75)(42,9){$Ep_{A_2}$}
\enddc}
\eqno(11)$$

and with every diagram

$$\vcenter
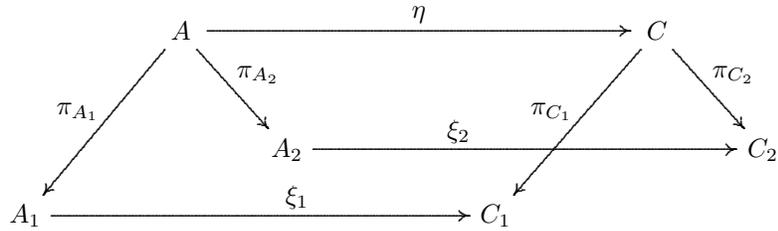
{
\begindc{0}[1]
\obj(1,50){$A_1$}
\obj(100,75){$A_2$}
\obj(60,120){$A$}
\mor(60,120)(100,75){$\pi_{A_2}$}
\mor(60,120)(1,50){$\pi_{A_1}$}[\atright,\solidarrow]
\obj(179,50){$C_1$}
\obj(280,75){$C_2$}
\obj(240,120){$C$}
\mor(240,120)(280,75){$\pi_{C_2}$}
\mor(240,120)(179,50){$\pi_{C_1}$}[\atright,\solidarrow]
\mor(1,50)(179,50){\phantom{ddddd}$\xi_1$}
\mor(100,75)(280,75){$\xi_2$\phantom{ddddddddd}}
\mor(60,120)(240,120){$\eta$}
\enddc}
$$

the diagram

$$\vcenter
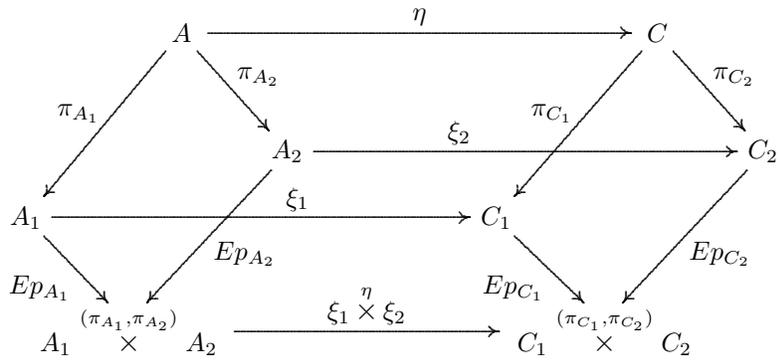
{
\begindc{0}[1]
\obj(40,7){$A_1\pout{(\pi_{A_1},\pi_{A_2})}A_2$}
\obj(1,50){$A_1$}
\obj(100,75){$A_2$}
\obj(60,120){$A$}
\mor(1,50)(40,10){$Ep_{A_1}$}[\atright,\solidarrow]
\mor(60,120)(100,75){$\pi_{A_2}$}
\mor(60,120)(1,50){$\pi_{A_1}$}[\atright,\solidarrow]
\mor(100,75)(40,10){$Ep_{A_2}$}
\obj(220,7){$C_1\pout{(\pi_{C_1},\pi_{C_2})}C_2$}
\obj(179,50){$C_1$}
\obj(280,75){$C_2$}
\obj(240,120){$C$}
\mor(240,120)(280,75){$\pi_{C_2}$}
\mor(240,120)(179,50){$\pi_{C_1}$}[\atright,\solidarrow]
\mor(179,50)(220,10){$Ep_{C_1}$}[\atright,\solidarrow]
\mor(280,75)(220,10){$Ep_{C_2}$}
\mor(70,7)(190,7){$\xi_1\pout{\eta}\xi_2$}
\mor(1,50)(179,50){\phantom{ddddd}$\xi_1$}
\mor(100,75)(280,75){$\xi_2$\phantom{ddddddddd}}
\mor(60,120)(240,120){$\eta$}
\enddc}
$$

\proclaim Proposition.
Diagram (11)  fulfills the universal property of pushouts in both categories $\k A$ and $\k C$.

{\it Proof.} \enspace
Suppose we have a diagram 

$$\vcenter{
\begindc{0}[1]
\obj(100,10){$D$}
\obj(10,10){$A_1$}
\obj(100,50){$A_2$}
\obj(10,50){$A$}
\mor(10,10)(100,10){$\sigma_1$}[\atright,\solidarrow]
\mor(10,50)(100,50){$\pi_{A_2}$}
\mor(10,50)(10,10){$\pi_{A_1}$}[\atright,\solidarrow]
\mor(100,50)(100,10){$\sigma_2$}
\enddc}
$$
where $\pi_{A_1}$ and $\pi_{A_2}$ are epimorphisms.
We first observe that
$$\ker\pi=\ker\pi_{A_1}+\ker\pi_{A_2}\subset\ker(\sigma_1\circ\pi_{A_1})=\ker(\sigma_2\circ\pi_{A_2})\,.$$
Then, from the Factorization Theorem, we conclude that there is a unique local algebra morphism $\tau$ for which the diagram 

$$\vcenter{
\begindc{0}[1]
\obj(10,10){$A/(\ker\pi_{A_1}+\ker\pi_{A_2})$}
\obj(120,50){$D$}
\obj(10,50){$A$}
\mor(10,50)(120,50){$\sigma_1\circ\pi_{A_1}=\sigma_2\circ\pi_{A_2}$}
\mor(10,50)(10,10){$\pi$}[\atright,\solidarrow]
\mor(15,15)(120,50){$\tau$}[\atright,\solidarrow]
\enddc}
$$

is commutative, i.e., $\tau\circ\pi=\sigma_1\circ\pi_{A_1}=\sigma_2\circ\pi_{A_2}$, or
$$\tau\circ Ep_{A_i}\circ\pi_{A_i}=\sigma_i\circ\pi_{A_i} \quad\quad (i=1,2)\,,$$
or, since $\pi_{A_i}$ are surjective, 
$$\tau\circ Ep_{A_i}= \sigma_i \quad\quad (i=1,2)\,.$$
This completes the proof.
\hfill$\bull$

We will call $\pout{(\cdot,\cdot)}$ the {\it pushout functor} in $\k A$ and we will refer to $A_1\pout{(\pi_{A_1},\pi_{A_2})}A_2$ as the pushout of $\pi_{A_1}$ and $\pi_{A_2}$. 

It is possible to consider the pushout of an arbitrary number of local algebra epimorphisms following the usual categorial method.
The associative property is trivially fulfilled.

\vskip0.5\baselineskip
(iii)\enspace  A Special Case. \enspace
In this subsection we will apply the pullback and pushout functors to the natural epimorphisms between quotients of a given local algebra.
We will thus establish a connection between the pullback and the intersection of ideals, and between the pushout and the sum of ideals in the sense that, if $A$ is a local algebra, and $J_1, J_2$ are proper ideals of $A$, the diagram of natural epimorphisms

$$\vcenter{
\begindc{0}[13]
\obj(4,1){$A/(J_1 + J_2)$}
\obj(1,4){$A/J_1$}
\obj(9,6){$A/J_2$}
\obj(6,9){$A/(J_1\cap J_2)$}
\mor(1,4)(4,1){$\Pi_1$}[\atright,\solidarrow]
\mor(6,9)(9,6){$\pi_2$}
\mor(6,9)(1,4){$\pi_1$}[\atright,\solidarrow]
\mor(9,6)(4,1){$\Pi_2$}
\enddc}
$$
is isomorphic to the pullback--pushout diagram

$$\vcenter{
\begindc{0}[13]
\obj(8,1){$A/J_1\pout{(\pi_1,\pi_2)}A/J_2\simeq A/J_1\pout{(Pr_1,Pr_2)}A/J_2$}
\obj(1,4){$A/J_1$}
\obj(9,6){$A/J_2$}
\obj(10,9){$A/J_1\pback{(\Pi_1,\Pi_2)}A/J_2\simeq A/J_1\pback{(Ep_1,Ep_2)}A/J_2$}
\mor(1,4)(4,1){$Ep_1$}[\atright,\solidarrow]
\mor(6,9)(9,6){$Pr_2$}
\mor(6,9)(1,4){$Pr_1$}[\atright,\solidarrow]
\mor(9,6)(4,1){$Ep_2$}
\enddc}
$$

In fact, according to (10), we have
$$
\eqalign{A/J_1\pout{(\pi_1,\pi_2)}A/J_2 &=  {A/(J_1\cap J_2) \over \ker \pi_1 + \ker \pi_2}
= {A/(J_1\cap J_2) \over J_1/(J_1\cap J_2) + J_2/(J_1\cap J_2)} = {A/(J_1\cap J_2) \over (J_1+ J_2)/(J_1\cap J_2)} \cr
&\simeq A/(J_1 + J_2),\cr
}$$
where the last isomorphism follows from the third isomorphism theorem for rings. 

Moreover, according to (8), we have
$$A/J_1 \pback{(\Pi_1,\Pi_2)}A/J_2 = \R +
 \{({\u a}_1 + J_1 ,{\u a}_2 + J_2) \in I_A/J_1\times I_A/J_2\mid {\u a}_1 + (J_1+J_2) = {\u a}_2 + (J_1+J_2)\}$$
and, as  a consequence of the first isomorphism theorem for rings, the mapping
$$\eqalign{
A/(J_1\cap J_2) &\lra  A/J_1 \pback{(\Pi_1,\Pi_2)}A/J_2\cr
a+(J_1\cap J_2)&\mapsto (a+J_1, a+J_2)\cr
}$$
is a local algebra isomorphism.

\Sezione {Composition of Functors.}

We conclude our study of the category of local algebras by examining the behaviour under composition of the functors introduced above together with the tensor product functor, which was considered by A. Weil [15].
We just recall that for any local algebras $A_1$ and $A_2$, the tensor product $A_1\otimes A_2$ is a local algebra, whose maximal ideal is $I_{A_1\otimes A_2}=I_{A_1}\otimes A_2+A_1\otimes I_{A_2}.$

We will only take into consideration the cases involving $\pback{(\cdot,\cdot)}$ with $\otimes$ and $\prod{(\cdot,\cdot)}$ with $\otimes$.

Given a local algebra $A$ and a diagram

\centerline{
\begindc{0}[13]
\obj(4,1){$B$}
\obj(1,4){$A_1$}
\obj(9,6){$A_2$}
\mor(1,4)(4,1){$\varphi_{A_1}$}[\atright,\solidarrow]
\mor(9,6)(4,1){$\varphi_{A_2}$}
\enddc}

we can consider the diagram

\centerline{
\begindc{0}[13]
\obj(4,1){$A\otimes B$}
\obj(1,4){$A\otimes A_1$}
\obj(9,6){$A\otimes A_2$}
\mor(1,4)(4,1){$Id_A\otimes\varphi_{A_1}$}[\atright,\solidarrow]
\mor(9,6)(4,1){$Id_A\otimes\varphi_{A_2}$}
\enddc}

\proclaim Proposition.
The functor $\otimes$ is distributive with respect to the functor $\pback{(\cdot,\cdot)}$, in the sense that  we have a natural local algebra isomorphism
$$A\otimes A_1\pback{(Id_A\otimes\varphi_{A_1},Id_A\otimes\varphi_{A_2})}A\otimes A_2 \simeq
A\otimes \left(A_1 \pback{(\varphi_{A_1},\varphi_{A_2})}A_2 \right)\eqno(12)$$

{\it Proof.}\enspace
As we already noticed in 2.3(i), the left-hand side of (12) is (isomorphic to) a subalgebra of $(A\otimes A_1)\times (A\otimes A_2)$.

Let $(e_1,\dots,e_d)$ be a basis  in $A$. We have that, for all $a_{i1}\in A_1$ and $a_{i2}\in A_2$ $(i=1,\dots,d)$, 

$$\left( \sum_{i=1}^d\, e_i\otimes a_{i1} , \sum_{i=1}^d\, e_i\otimes a_{i2} \right) \in
A\otimes A_1\pback{(Id_A\otimes\varphi_{A_1},Id_A\otimes\varphi_{A_2})}A\otimes A_2$$
if and only if
$$\sum_{i=1}^d\, e_i\otimes \varphi_{A_1}(a_{i1}) =  \sum_{i=1}^d\, e_i\otimes  \varphi_{A_2}(a_{i2})$$
which is equivalent to
$$\varphi_{A_1}(a_{i1}) = \varphi_{A_2}(a_{i2}) \quad\quad \forall\, i=1,\dots,d$$ 
i.e., to
$$(a_{i1},a_{i2})\in A_1 \pback{(\varphi_{A_1},\varphi_{A_2})}A_2  \quad\quad \forall\, i=1,\dots,d$$
and, finally to
$$\sum_{i=1}^d\,e_i\otimes(a_{i1},a_{i2})\in A\otimes \left(A_1 \pback{(\varphi_{A_1},\varphi_{A_2})}A_2\right)\,.$$
It can be easily shown that the resulting local algebra isomorphism is independent of the choice of the basis in $A$.
\hfill$\bull$

\vskip0.5\baselineskip
Recall that $\times_{\R}$ is a special case of $\pback{(\cdot,\cdot)}$ and that $\prod{(\cdot,\cdot)}$ is obtained by means of  triple products. 
Therefore, as a trivial consequence of the above proposition, we have 

\proclaim Corollary.
The functor $\otimes$ is distributive with respect to both the functors $\times_{\R}$ and $\prod{(\cdot,\cdot)}$.
\hfill$\bull$

It is easy to check that functor $\times_{\R}$ is distributive with respect to functor $\pback{(\cdot,\cdot)}$. 

\Capitolo {Algebras of Truncated Polynomials.}

After studying the category of local algebras in general, the next natural step is to look for concrete examples.
The model of the most general local algebra is a finite dimensional quotient of an algebra of polynomials over a local algebra $A$, a so-called algebra of \lq generalized truncated polynomials'.
We will apply some of the most significant functors introduced in the previous section in order to obtain special algebras of this kind.
In particular we will find examples of local algebras which are isomorphic as $A$--modules but not as algebras.

\vskip0.5\baselineskip
Let $A$ be a $d$-dimensional local algebra.

Consider the algebra $A\lbrack x\rbrack$ of all polynomials with coefficients in $A$ in the indeterminate $x$.

Let $\II=\langle x\rangle$
be the ideal of $A\lbrack x\rbrack$ generated by $x$; then, for all non negative integers $k$, 
$\II^k=\langle x^k\rangle$.

The proof of the next proposition is reported for the only sake of completeness. 

\proclaim Proposition.
The set of all residue classes mod $x^{k+1}$ of polynomials in $A\lbrack x\rbrack$, i.e.,
$$\PP_k A\lbrack x\rbrack:=A\lbrack x\rbrack/\II^{k+1}\,,$$
is a local algebra.

{\it Proof.}\enspace
Since $\dim_{\R}(\PP_k A\lbrack x\rbrack)=(k+1)d$, condition (i) of Definition 1 is fulfilled.

Now consider a maximal ideal $J$ of $A\lbrack x\rbrack$ containing $\II^{k+1}$.  
Since $J$ is prime and $x^{k+1}\in J$, we have that $x\in J$.
Hence $\II\subset  J$, so that  $J=I+\II$ for some proper ideal $I$ of $A$.
Therefore $J\subset I_A+\II$, whence equality follows by maximality.
We have thus proved that $\PP_k A\lbrack x\rbrack$ has $I_k=(I_A+\II/\II ^{k+1})$ as its unique maximal ideal. 
It also holds that $\PP_k A\lbrack x\rbrack/I_k\simeq {\R}$.
Hence condition (ii) of Definition 1 is fulfilled, too.
\hfill$\bull$

\vskip0.5\baselineskip
$\PP_k A\lbrack x\rbrack$ can be identified with the set of all polynomials of degree at most $k$.

\vskip1\baselineskip

We can consider the following algebras:
$$\PP_k A\lbrack x_1,\dots, x_n\rbrack:=A\lbrack x_1,\dots, x_n\rbrack/\II^{k+1}\;, \eqno(13)$$
where $\II=\langle x_1,\dots, x_n\rangle$, and
$$\PP_{k_1,\dots,k_n} A\lbrack x_1,\dots, x_n\rbrack:=A\lbrack x_1,\dots, x_n\rbrack/\II^{(k_1+1,\dots,
k_n+1)}\eqno(14)$$
where
$$\II^{(k_1+1,\dots,k_n+1)}:=\langle x_1^{k_1+1},\dots,x_n^{k_n+1}\rangle\;.$$
By the same arguments as in Proposition 6, one can easily conclude that all the above algebras (13) and (14) are local.
We will call them {\it algebras of generalized truncated polynomials}.

We now apply the categorial constructions described in the previous section to algebras of generalized truncated polynomials. 
In all the cases considered, the resulting local algebras will still be of the same kind. 

\vskip0.5\baselineskip
Note that we have the following inclusions:

$$\II^k\subset \II^{k'}\iff k\ge k'$$
$$\II^{(k_1,\dots,k_n)}\subset \II^{(k'_1,\dots,k'_n)}\iff k_i\ge k'_i\qquad\hbox{for all }i=1,\dots,n $$
$$\II^k\subset \II^{(k_1,\dots,k_n)}\iff k\geq k_1+\cdots +k_n-n+1$$
$$ \II^{(k_1,\dots,k_n)}\subset\II^k\iff k_i\geq k\qquad\hbox{for all }i=1,\dots,n$$

In each of the above cases we can consider the corresponding quotients of the local algebras defined in (13) and (14):
$$\PP_{k-1} A\lbrack x_1,\dots, x_n\rbrack / (\II^{k'} /  \II^k) \simeq \PP_{k'-1} A\lbrack x_1,\dots, x_n\rbrack $$

$$\PP_{(k_1-1,\dots,k_n-1)} A\lbrack x_1,\dots, x_n\rbrack / \left(\II^{(k'_1,\dots,k'_n)} /  \II^{(k_1,\dots,k_n)}\right) 
\simeq \PP_{(k'_1-1,\dots,k'_n-1)} A\lbrack x_1,\dots, x_n\rbrack$$

$$\PP_{k-1} A\lbrack x_1,\dots, x_n\rbrack /  \left(\II^{(k_1,\dots,k_n)}/ \II^k \right)
 \simeq \PP_{(k_1-1,\dots,k_n-1)} A\lbrack x_1,\dots, x_n\rbrack  $$
 
 $$\PP_{(k_1-1,\dots,k_n-1)} A\lbrack x_1,\dots, x_n\rbrack / \left(\II^k /  \II^{(k_1,\dots,k_n)}\right) 
\simeq \PP_{k-1} A\lbrack x_1,\dots, x_n\rbrack$$

All the above isomorphisms are a consequence of the third isomorphism theorem for rings.

\vskip0.5\baselineskip
We also have the natural isomorphisms

$$\bigotimes_{i=1}^n\PP_{k_i}A[x_i]\simeq\PP_{k_1,\dots,k_n} A[x_1,\dots,x_n]$$

\vskip0.5\baselineskip

Let $m\in\Bbb{N}^*$. 

We will show that the free $A$--module $A^m$ of rank $m$ can be given different local algebra structures.

For any $r,s,t\in\Bbb{N}^*$ such that $r+t-s+1=m$ and $s\leq\min(r,t)$, consider the diagram

\vskip0.5\baselineskip
$$\vcenter{
\begindc{0}[3]
\obj(17,10){$\PP_tA\lbrack y\rbrack$}
\obj(42,10){$\PP_sA\lbrack z\rbrack$}
\obj(42,25){$\PP_rA\lbrack x\rbrack$}
\mor(20,10)(38,10){$\tau^s_t$}[\atright,\solidarrow]
\mor(40,25)(40,10){$\tau^s_r$}
\enddc}\eqno(15)$$

\vskip0.5\baselineskip

where $\tau^s_r\big(\sum_{i=0}^ra_ix^i\big)=\sum_{i=0}^sa_iz^i$, and $\tau^s_t$ is defined similarly (with a slight abuse of notation), and consider the corresponding (relative) product 
 $$\PP_rA\lbrack x\rbrack \prod{(\tau^s_r,\tau^s_t)}\PP_tA\lbrack y\rbrack =
 \PP_sA\lbrack z\rbrack \times_{\R}{\rm Ker}\,\tau^s_r \times_{\R}{\rm  Ker}\,\tau^s_t\,.\eqno(16)$$
We introduce the free $A$--module of rank $m$
$$C^s_{r,t}=\left\{ \left( a_0+\sum_{i=1}^sa_iz^i \right) + \sum_{i=s+1}^rb_ix^i+\sum_{i=s+1}^tc_iy^i\right\},$$
which can be identified in an obvious way with a subset of
$$\PP_{r,t,s}A\lbrack x,y,z\rbrack/\langle xy,xz,yz\rangle\,.$$
It is easy to see that $C^s_{r,t}$, endowed with the  induced local algebra structure, is isomorphic to (16).

\vskip0.5\baselineskip

We now apply the pullback functor to diagram (15) and consider the corresponding local algebra
$$\PP_rA\lbrack x\rbrack \pback{(\tau^s_r,\tau^s_t)}\PP_tA\lbrack y\rbrack =
{\R}+\left\{(\u p_1,\u p_2)\in I_{\PP_rA\lbrack x\rbrack }\times I_{\PP_tA\lbrack y\rbrack}
\mid \tau^s_r(\u p_1) = \tau^s_t(\u p_2)\right\} .\eqno(17)$$
We introduce the free $A$--module of rank $m$
$$B^s_{r,t}=\left\{a_0+\sum_{i=1}^sa_i(x^i+y^i)+\sum_{i=s+1}^rb_ix^i+\sum_{i=s+1}^tc_iy^i\right\},$$
which can be identified, in an obvious way, with a subset of
$$\PP_{r,t}A\lbrack x,y\rbrack/\langle xy\rangle\,.$$
It is easy to see that $B^s_{r,t}$, endowed with the  induced local algebra structure, is isomorphic to (17).

\vskip0.5\baselineskip
In both of the above constructions, different decompositions of integer $m$ produce different $A$--modules (of the same rank $m$), which, in general, are pairwise non-isomorphic as local algebras. This yields different local algebra structures on $A^m$.

\Capitolo {Points Proches.}

We start this section by recalling the construction of the functor associated with a local algebra, introduced by A.Weil [15].

\Sezione {Functors associated with Local Algebras.}

Let $A$ be a local algebra.
For any differential manifold $P$ consider the set $AP$ of all morphisms $u: C^{\infty}(P)\to A$ which preserve the identity.
The elements in $AP$  are called {\it points proches} associated with the local algebra $A$.

Let $u\in AP$. 
There is just one point $p\in P$ such that the composition
$$0_A\circ u : C^{\infty}(P)\lra {\R}$$
is the evaluation mapping $ev_p$ at $p$ (cf., e.g., [7], p.296 ).

This defines a surjective mapping ({\it target}) $\alpha^T_P : AP\lra P$ such that
$$(\alpha^T_P)^{-1}(p)=A_pP=\big\{u\in AP \vert 0_A\circ u=ev_p\big\}\, .$$

Now consider a differentiable mapping $\varphi: P\lra Q$ and define
$$A\varphi : AP\lra AQ$$
by setting
$$A\varphi(u):C^{\infty}(Q)\lra \R  : g\mapsto u(g \circ \varphi).$$

For each $p\in P$, this mapping induces by restriction  a mapping
$$A_p\varphi : A_pP\lra A_{\varphi(p)}Q\,.$$

\proclaim Proposition.
\endgraf
(i)\enspace To any local algebra $A$ corresponds the covariant functor 
$$\vcenter{
\begindc{0}[1]
\obj(1,1){$I$}
\obj(1,50){$A$}
\mor(1,50)(1,1){$\alpha^T$}[\atright,\solidarrow]
\enddc}
$$
which associates with any differential manifold $P$ the mapping $\alpha^T_P : AP\lra P$
and with any differentiable mapping $\varphi: P\to Q$ the diagram
$$\vcenter{
\begindc{0}[1]
\obj(1,1){$P$}
\obj(1,50){$AP$}
\obj(100,1){$Q$}
\obj(100,50){$AQ$}
\mor(1,50)(1,1){$\alpha_P^T$}[\atright,\solidarrow]
\mor(100,50)(100,1){$\alpha_Q^T$}
\mor(1,1)(100,1){$\varphi$}
\mor(1,50)(100,50){$A\varphi$}
\enddc}
$$
(ii) \enspace To a morphism of local algebras $\kappa:A\lra B$ corresponds the functor
$$\vcenter{
\begindc{0}[1]
\obj(1,1){$I$}
\obj(1,50){$A$}
\obj(100,1){$I$}
\obj(100,50){$B$}
\mor(1,50)(1,1){$\alpha^T$}[\atright,\solidarrow]
\mor(100,50)(100,1){$\beta^T$}
\mor(1,1)(100,1){$$}[\atright,\solidline]
\mor(1,4)(100,4){$$}[\atright,\solidline]
\mor(1,50)(100,50){$\kappa$}
\enddc}
$$
which associates with a manifold $P$ the diagram
$$\vcenter{
\begindc{0}[1]
\obj(1,1){$P$}
\obj(1,50){$AP$}
\obj(100,1){$P$}
\obj(100,50){$BP$}
\mor(1,50)(1,1){$\alpha^T_P$}[\atright,\solidarrow]
\mor(100,50)(100,1){$\beta^T_P$}
\mor(1,1)(100,1){$$}[\atright,\solidline]
\mor(1,4)(100,4){$$}[\atright,\solidline]
\mor(1,50)(100,50){$\kappa_P$}
\enddc}
$$
where
$$\kappa_P(u)=\kappa\circ u\,,$$
and with a morphism $\varphi:P\to Q$ the diagram

\centerline{
\begindc{0}[1]
\obj(40,10){$P$}
\obj(170,10){$Q$}
\obj(10,40){$P$}
\obj(140,40){$Q$}
\obj(90,60){$BP$}
\obj(220,60){$BQ$}
\obj(60,90){$AP$}
\obj(190,90){$AQ$}
\mor(40,10)(170,10){$\varphi$}
\mor(10,40)(140,40){$\phantom{bbbbbbbb}\varphi$}
\mor(90,60)(220,60){$B\varphi\phantom{bbbbbbbbbb}$}
\mor(60,90)(190,90){$A\varphi$}
\mor(10,40)(40,10){$$}[\atright,\solidline]
\mor(14,40)(44,10){$$}[\atright,\solidline]
\mor(140,40)(170,10){$$}[\atright,\solidline]
\mor(144,40)(174,10){$$}[\atright,\solidline]
\mor(60,90)(90,60){$\kappa_P$}
\mor(190,90)(220,60){$\kappa_Q$}
\mor(60,90)(10,40){$\alpha^T_P$}[\atright,\solidarrow]
\mor(90,60)(40,10){$\beta^T_P$}
\mor(190,90)(140,40){$\alpha^T_Q$}[\atright,\solidarrow]
\mor(220,60)(170,10){$\beta^T_Q$}
\enddc}
\hfill$\bull$

In the first part of the paper we introduced some basic categorial constructions in order  to obtain new local algebras from old ones. In view of the above proposition, to each of these algebras corresponds a functor.
Weil noticed that to the tensor product of local algebras corresponds the composition of the functors associated with the factors. 
Here we study what corresponds to the pullback of local algebras.

Let

\vskip0.5\baselineskip
\centerline{
\begindc{0}[1]
\obj(90,60){$A_1$}
\obj(220,60){$B$}
\obj(60,100){$A_1\pback{(\varphi_{A_1},\varphi_{A_2})} A_2$}
\obj(190,100){$A_2$}
\mor(90,60)(220,60){$\varphi_{A_1}$}[\atright,\solidarrow]
\mor(90,100)(190,100){$Pr_{A_2}$}
\mor(60,95)(90,60){$Pr_{A_1}$}[\atright,\solidarrow]
\mor(190,95)(220,60){$\varphi_{A_2}$}
\enddc}
\vskip0.5\baselineskip

be the pullback diagram of two local algebra morphisms $\varphi_{A_1}:A_1\to B$ and $\varphi_{A_2}:A_2\to B$.

Then, for any manifold $P$, we have the diagram

\vskip1\baselineskip
\centerline{
\begindc{0}[1]
\obj(40,10){$P$}
\obj(170,10){$P$}
\obj(10,40){$P$}
\obj(140,40){$P$}
\obj(90,60){$A_1P$}
\obj(220,60){$BP$}
\obj(60,100){$(A_1\pback{(\varphi_{A_1},\varphi_{A_2})} A_2)P$}
\obj(190,100){$A_2P$}
\mor(40,10)(170,10){$$}[\atright,\solidline]
\mor(40,13)(170,13){$$}[\atright,\solidline]
\mor(10,40)(140,40){$$}[\atright,\solidline]
\mor(10,43)(140,43){$$}[\atright,\solidline]
\mor(90,60)(220,60){$\varphi_{A_1\,P}\phantom{aaaaaaaaa}$}[\atright,\solidarrow]
\mor(90,95)(190,95){$Pr_{A_2P}$}
\mor(10,40)(40,10){$$}[\atright,\solidline]
\mor(14,40)(44,10){$$}[\atright,\solidline]
\mor(140,40)(170,10){$$}[\atright,\solidline]
\mor(144,40)(174,10){$$}[\atright,\solidline]
\mor(60,95)(90,60){$Pr_{A_1P}$}
\mor(190,95)(220,60){$\varphi_{A_2P}$}
\mor(60,95)(10,40){${(\alpha_1\pback{(\varphi_{A_1},\varphi_{A_2})} \alpha_2)^T_P}$}[\atright,\solidarrow]
\mor(90,60)(40,10){$\alpha^T_{1P}$}
\mor(190,95)(140,40){$\alpha^T_{2P}$}[\atright,\solidarrow]
\mor(220,60)(170,10){$\beta^T_P$}
\enddc}
\vskip1\baselineskip

On the other hand, consider the pullback diagram of $\varphi_{A_1\,P}$ and $\varphi_{A_2\,P}$

\vskip0.5\baselineskip
\centerline{
\begindc{0}[1]
\obj(90,60){$A_1P$}
\obj(220,60){$BP$}
\obj(60,100){$A_1P\times_{BP}A_2P $}
\obj(190,100){$A_2P$}
\mor(92,60)(220,60){$\varphi_{A_1P}$}[\atright,\solidarrow]
\mor(90,100)(188,100){$Pr_{A_2P}$}
\mor(60,95)(90,60){$Pr_{A_1P}$}[\atright,\solidarrow]
\mor(190,95)(220,60){$\varphi_{A_2P}$}
\enddc}
\vskip0.5\baselineskip
where
$$A_1P\times_{BP}A_2P = \{ (u_1,u_2)\in A_1P\times A_2P \vert \varphi_{A_1\,P}(u_1) = \varphi_{A_2\,P}(u_2) \}.$$

If we put
$$\tau = \varphi_{A_1P} \circ Pr_{A_1P} = \varphi_{A_2P} \circ Pr_{A_2P}\,,$$
it is easy to prove that in the following diagram

\vskip1\baselineskip
\centerline{
\begindc{0}[1]
\obj(30,170){$P$}
\obj(200,170){$P$}
\obj(30,230){$(A_1\pback{(\varphi_1,\varphi_2)} A_2)P$}
\obj(200,230){$A_1P\times_{BP} A_2P$}
\mor(55,230)(175,230){$(Pr_{A_1P},Pr_{A_2P})$}[\atleft,\solidarrow ]
\mor(30,171)(200,171){$$}[\atright,\solidline]
\mor(30,168)(200,168){$$}[\atright,\solidline]
\mor(30,230)(30,170){$(\alpha_1\pback{(\varphi_{A_1},\varphi_{A_2})} \alpha_2)^T_P$}[\atright,\solidarrow ]
\mor(200,230)(200,170){$\beta^T_P\circ\tau$}
\enddc}
\vskip0.5\baselineskip

the upper arrow is a bijection. 

This allows us to conclude that the pullback commutes with the operation of associating with a local algebra the corresponding points proches. 

Since  the product of local algebras is a special case of pullback, we have the diagram
\vskip0.5\baselineskip
\centerline{
\begindc{0}[1]
\obj(30,170){$P$}
\obj(200,170){$P$}
\obj(30,230){$(A_1\times_{\R}A_2)P$}
\obj(200,230){$A_1P\times_{P} A_2P$}
\mor(55,230)(175,230){$(Pr_{A_1P},Pr_{A_2P})$}[\atleft,\solidarrow ]
\mor(30,175)(200,173){$$}[\atright,\solidline]
\mor(30,170)(200,170){$$}[\atright,\solidline]
\mor(30,230)(30,170){$(\alpha_1\times_{\R}\alpha_2)^T_P$}[\atright,\solidarrow]
\mor(200,230)(200,170){$\tau$}
\enddc}
\vskip0.5\baselineskip
and the upper arrow is still  bijective. 
In this case,
$$A_1P\times_{P} A_2P = \{ (u_1,u_2)\in A_1P\times A_2P \vert \alpha^T_{1\,P}(u_1) = \alpha^T_{2\,P}(u_2) \}\,,$$
is the fibre product of the target mappings.

We conclude that the product of local algebras gives rise to the fibre product of the corresponding points proches.

\Sezione  {The Local Character of Points Proches.}

Let $A$ be a local algebra and $P$ a differential manifold.

\proclaim Proposition.
Let $u\in A_pP$. Then $u(f)=0$, for any function $f$ vanishing in some neighbourhood of $p$. 

{\it Proof.}\enspace 
Let $f$ be a function vanishing in the neighbourhood $\k U$ of $p$. 
If $\beta$ is a bump function at $p$ with support in $\k U$, we have
$$f = (1-\beta)f$$
both on $\k U$, where $f$ vanishes, and on $P-\k U$, where $\beta$ vanishes. 
Therefore
$$u(f)=u(1-\beta)u(f)\, .\eqno(18)$$
Put
$$\eqalignno{
&u(f)=f(p)+\u a= \u a \cr
&u(1-\beta)=(1-\beta)(p)+\u b = \u b\, .\cr
}$$
Owing to (18) we have
$$\u a=\u b\, \u a\iff (1-\u b)\u a=0\,,$$
but, since $\u b\in I_A$ and $1\notin I_A$, it must be $1-\u b\notin I_A$. 
This implies that $1-\u b$ is invertible.
Hence $u(f) = \u a = 0$.
\hfill$\bull$

\vskip0.5\baselineskip
Proposition 8 reveals the local character of points proches.
It allows us to prove that the functor associated with a local algebra $A$ preserves open inclusions.
Let us apply this functor to the inclusion mapping
$$\iota_{\k U}  : \k U\lra P$$
of the open subset $\k U$ into $P$.
Then for each $p\in \k U$, we have the mapping
$$A_p\iota_{\k U}  : A_p\k U\lra A_pP$$
given by
$$A_p\iota_{\k U}(u)(f):=u(f\vert _{\k U})\,.$$
On the other hand, the mapping
$$\eta : A_pP \lra A_p\k U $$
given by
$$ \eta(u)(g) := u(f)$$
for any  differentiable prolongation $f$ of $g$ in $P$, is well defined, in view of the local character of $u$,  and is the inverse of $A_p\iota_{\k U}$. 
Hence $A_p\iota_{\k U}$ is a bijection.

If we identify $A_p{\k U}$ with $A_pP$
and consequently regard $A\iota_{\k U}$ as the inclusion $\iota_{A\k U}$ of $A\k U$ in $AP$, we have the following diagram:

$$\vcenter{
\begindc{0}[3]
\obj(10,10){$\k U$}
\obj(40,10){$P$}
\obj(10,25){$A\k U$}
\obj(40,25){$AP$}
\mor(10,24)(40,24){$\iota_{A\k U}$}
\mor(10,25)(10,10){$\alpha^T_{\k U}$}[\atright,\solidarrow]
\mor(40,25)(40,10){$\alpha^T_P$}
\mor(10,9)(40,9){$\iota_{\k U}$}
\enddc}\eqno(19)
$$

\vskip0.5\baselineskip
This result will be used for introducing the differential fibration structure on $AP$.

\Sezione {Points Proches in Affine Spaces.}

In this section we study the special case where the functor associated with a local algebra is applied to an affine space. 
We will prove that the result is a fibre bundle whose total space has an affine structure.
In addition, we will show that the functor transforms every differentiable mapping between affine spaces into a differentiable morphism of trivial bundles.

\proclaim Proposition.
Let $(M,V)$ be an affine space.
Then, for each $x_0\in M$, we have a natural bijection
$$A_{x_0}M \lra \{x_0\} \times I_A\otimes V\,.$$

{\it Proof.}\enspace
We denote by $K^s(M,x_0)$ the linear space of homogeneous polynomials in $M$, at $x_0$, of degree $s$ (cf. Definition 13 in the Appendix).
There is a natural isomorphism
$$K^1(M,x_0)\simeq V^*$$
which allows us to identify $I_A\otimes V$ with ${\rm Hom}_{\R}(K^1(M,x_0),I_A)$.
We will prove that there is a natural bijection
$$A_{x_0}M \lra  \{x_0\} \times{\rm Hom}_{\R}\big(K^1(M,x_0),I_A \big).$$
Let $(e_1,\dots,e_m)$ be a basis of $V$, $(e^1,\dots,e^m)$ its dual basis and 
$$x^i : M\lra{\R} : x\mapsto e^i(x-x_0) \eqno(20)$$
the corresponding Cartesian coordinate functions at $x_0$, for $i=1,\dots,m$.  

Let us now consider the mapping
$$\eqalign{
A_{x_0}M &\lra  \{x_0\} \times{\rm Hom}_{\R}\big(K^1(M,x_0),I_A \big)\cr
u & \mapsto \left( x_0, u\vert_{K^1(M,x_0)} \right) \cr
}$$

We will prove that it admits an inverse, i.e., for each $\bar u\in {\rm Hom}_{\R}\big(K^1(M,x_0),I_A\big)$ there is a unique $u\in A_{x_0}M$ with $u\vert_{K^1(M,x_0)} = \bar u$. 
We are going to construct $u$ starting from its action on polynomials.
Let $h\in K^s(M,x_0)$. 
In view of Proposition 14 in the Appendix, we have
$$\eqalign{
h(x) & = {1\over s!} \Delta^s h(x_0 ; x-x_0)\cr
&= {1\over s!}\Delta^s h(x_0 ; x^ie_i)\cr
&= {1\over s!}\delta^s h(x_0 ; e_{i_1},\dots ,e_{i_s})x^{i_1}\cdots x^{i_s}\cr
}$$
where $\Delta^s h$ and $\delta^s h$ are the $s$-th unidirectional and (multidirectional) polarizations of $h$, respectively (cf. Definition 12 in the Appendix) and we have adopted the Einstein convention over repeated indices.
We define
$$ u(h) :=  {1\over s!}\delta^s h(x_0 ; e_{i_1},\dots ,e_{i_s})\bar u(x^{i_1})\cdots \bar u(x^{i_s})\,.$$
Easy linear algebraic arguments  show that this definition is independent of the choice of the basis.
It follows that 
$$u(h)=0 \hbox{ if } s >\ell $$
where $\ell$ is the height of $A$, and that $u$ preserves the product of polynomials.

We extend the definition of $u$ from polynomials to arbitrary smooth functions $f\in C^{\infty}(M)$ by considering the Taylor polynomial of $f$ of order $\ell$ at $x_0$ (cf. Note 18 in the Appendix) and setting
$$u(f) := f(x_0) + \sum_{i=1}^{\ell}\;{1\over i!} u\big( D^if(x_0 ; \cdot - x_0) \big) \,. \eqno(21)$$ 
We recall that each $D^if(x_0 ; \cdot - x_0)$ is a homogeneous polynomial at $x_0$ of degree $i$ (cf. Proposition 16 in the Appendix). 

The above formula defines $u$ as a linear mapping.
We show that $u$ preserves products, too. 
In fact, 
$$ u(f)u(g) = \left[ f(x_0) + \sum_{r=1}^{\ell}\;{1\over r!} u\Big( D^rf(x_0 ; \cdot - x_0) \Big)\right]
       \left[ g(x_0) + \sum_{s=1}^{\ell}\;{1\over s!} u\Big( D^sg(x_0 ; \cdot - x_0) \Big)\right].$$
Note that
$$u\Big(D^r(f)(x_0;\cdot - x_0)\Big) u\Big(D^s(g)(x_0;\cdot-x_0) \Big)= 0\quad {\rm if} \quad r+s>\ell.$$
Hence, by the Leibniz rule we deduce that
$$\eqalign{
u(f)u(g) 
&=f(x_0)g(x_0) + \sum_{r+s=1}^{\ell} {1\over r!s!}u\Big(D^r(f)(x_0;\cdot - x_0)\Big) u\Big(D^s(g)(x_0;\cdot-x_0) \Big) \cr 
&= f(x_0)g(x_0) + \sum_{i=1}^{\ell}\;{1\over i!}
     u\left( \sum_{r+s=i} {i!\over r!s!}D^r(f)(x_0;\cdot - x_0)D^s(g)(x_0;\cdot-x_0) \right)\cr 
&= fg(x_0) + \sum_{i=1}^{\ell} {1\over i!}\; u\big( D^i(fg)(x_0 ; \cdot - x_0) \big) \cr
&=u(fg). \cr
}$$

We have constructed a $u\in A_{x_0}M$ such that $u\vert_{K^1(M,x_0)} = \bar u$. 
It is the only element in $ A_{x_0}M$ fulfilling this property because, on the one hand, the action of an algebra morphism on polynomials is uniquely determined by its action on $K^1(M,x_0)$ and, on the other hand, each element in $ A_{x_0}M$  vanishes on all  homogeneous polynomials of degree $\ell +1$ at $x_0$; in particular it vanishes on the Lagrange remainders of order $\ell$ at $x_0$ (cf. Note 18 in the Appendix). 
As a consequence, its action on a function must be the one defined in (21).
\hfill$\bull$

\vskip0.5\baselineskip
According to Proposition 9 we have a natural bijection
$$AM \lra M\times I_A\otimes V \eqno(22)$$
which induces on  $AM$ a structure of $m(d-1)$--dimensional affine space, $d$ being the dimension of $A$.

By virtue of (22) we will identify $\alpha^T_M: AM\to M$ with
$$Pr_M : M\times I_A\otimes V \lra M$$

\proclaim Proposition.
Let $(M,V)$ and $(N,W)$ be affine spaces, and $\varphi : M\lra N$ a differentiable mapping. Then $A\varphi : AM\lra AN$ is a differentiable mapping.

{\it Proof.}\enspace
Let $(e)=(e_1,\dots,e_m)$ and $(\varepsilon)=(\varepsilon_1,\dots\varepsilon_n)$ be  bases of $V$ and $W$ respectively. 
We can identify $AM$ with
$$ M \times \big( \underbrace{I_A\times\dots\times I_A}_{\scriptstyle{m\,{\rm times}}}\big)$$
via the affine isomorphism
$$AM  \lra M \times \big( \underbrace{I_A\times\dots\times I_A}_{\scriptstyle{m\,{\rm times}}}\big)
\,:\, u  \mapsto \big(x, u(x^1),\dots,u(x^m)\big),$$
where $x=\alpha^T_M(u)$ and, $(x^1,\dots,x^m)$ are the Cartesian coordinate functions at $x$ relative to $(e)$, defined as in (20).
Similarly we identify $AN$ with $N\times \big( \underbrace{I_A\times\dots\times I_A}_{\scriptstyle{n\,{\rm times}}}\big)$.

Now consider $A\varphi : AM\to AN$ and recall that, for each $g\in C^{\infty}(N)$,
$$A\varphi(u)(g) = u(g\circ \varphi).$$
From the above identifications it follows that
$$u= \big(x, u(x^1),\dots,u(x^m)\big)$$
and
$$A\varphi(u)= \big(\varphi(x), A\varphi(u)(y^1),\dots,A\varphi(u)(y^n)\big),$$
where $(y^1,\dots,y^n)$ are the Cartesian coordinate functions at $\varphi(x)$ relative to $(\varepsilon)$.
Moreover, for all $j=1,\dots,n$ we have:
$$\eqalign{
A\varphi(u)(y^j) &= u(y^j \circ \varphi) \cr
    &= u(\varphi^j)\cr
    &= u \left[ \varphi^j(x) + \sum_{k=1}^{\ell} {1\over k!} D^k \varphi^j (x; \cdot - x) \right]\cr
    &= \varphi^j(x) + u \big( h^j [x^1,\dots,x^m] \big) \cr
}$$
where $h^j$ is a formal polynomial in $m$ indeterminates, whose coefficients are differentiable functions of $x$.
Hence
$$A\varphi(u)(y^j) = \varphi^j(x) + h^j [u(x^1),\dots, u(x^m)]$$
which shows that $A\varphi(u)(y^j)$ is differentiable as a function of $x$ and $u(x^i)$.
\hfill$\bull$

\vskip0.5\baselineskip
The above arguments imply that

\centerline{
\begindc{0}[3]
\obj(10,10){$M$}
\obj(40,10){$N$}
\obj(10,25){$AM$}
\obj(40,25){$AN$}
\mor(10,25)(40,25){$A\varphi$}[\atleft,\solidarrow]
\mor(10,25)(10,10){$\alpha^T_M$}[\atright,\solidarrow]
\mor(40,25)(40,10){$\alpha^T_N$}
\mor(10,10)(40,10){$\varphi$}[\atright,\solidarrow]
\enddc}

is a differentiable morphism of trivial bundles.

If $\varphi$ is a diffeomorphism, the diagram is an isomorphism.

In view of the local character of the points proches, all the above results are still true if the affine spaces are replaced by open subsets.

\Sezione {The differential structure.}

We now show that the functor associated with a local algebra transforms differential manifolds into differential fibre bundles.

\proclaim Proposition.
Let $P$  be a differential manifold modelled on an affine space $(M,V)$. 
Then $AP$ is a differential manifold modelled on the affine space $(AM,I_A \otimes V)$.

{\it Proof.} \enspace 
Let $\xi : \k U\subset P\to\xi(\k U)\subset M$ be an admissible chart of $P$.
We have the bijection
$$A\xi : A\k U \subset AP\lra A\xi(\k U)\subset AM\,.$$
where the inclusions are due to (19).
Now let $\xi' : \k U'\subset P\to\xi'(\k U')\subset M$ be another admissible chart on $P$ such that $\k U\cap\k U'\ne\emptyset$, then
$$\xi '\circ\xi ^{-1} : \xi (\k U\cap\k U') \lra \xi' (\k U\cap\k U')$$
is a diffeomorphism between open subsets of $M$.
By Proposition 11 
$$A(\xi '\circ\xi ^{-1}) = A\xi '\circ (A\xi )^{-1}: A\xi (\k U\cap\k U') \lra A\xi' (\k U\cap\k U')$$
is a diffeomorphism between open subsets of $AM$.
This guarantees that the functor associated with local algebra $A$, when applied to a differential atlas of $P$, produces a differential atlas of $AP$.
\hfill$\bull$

\vskip0.5\baselineskip
An immediate consequence of the above proposition is that the functor associated with a local algebra $A$, applied to chart $\xi$ of $P$, gives rise to a local trivialization of $AP$:
\vskip0.5\baselineskip

\centerline{
\begindc{0}[3]
\obj(10,10){$\k U$}
\obj(35,10){$\xi(\k U)$}
\obj(10,25){$A{\k U}$}
\obj(35,25){$\xi(\k U)\times I_A\otimes V$}
\mor(10,25)(26,25){$A\xi$}
\mor(10,10)(34,10){$\xi$}[\atright,\solidarrow]
\mor(10,25)(10,10){$\alpha^T_{\k U}$}[\atright,\solidarrow]
\mor(35,25)(35,10){$\alpha^T_{\xi(\k U})$}
\enddc}

so that we can conclude that
$$\alpha^T_P : AP \lra P$$ 
is a differentiable fibre bundle.

\Sezione  {Appendix.} We report from [1] some notions and results on real functions defined on affine spaces.

\proclaim Definition.	
Let $(Q,V)$ be a real affine space.											
The $0$-{\it th polarization} of a function $f : Q \to \R$ is the function itself.  For $n > 0$, the $n$-{\it
th polarization} of $f$ is the function
$$\delta^n f \colon Q \times V^n \to \R $$
defined by
$$
\delta^n f(q; v_1, \dots, v_n) = (-1)^n  \sum_{m=0}^n (-1)^m \!\!\!
\sum_{1\le i_1<\dots<i_m\le n} f(q + v_{i_1}  +\dots + v_{i_m}). 
$$ 
The term in the above sum corresponding to $m=0$ is set equal to $f(q)$.
\hfill$\bull$

We introduce the $n$-{\it th unidirectional polarization} of a function, defined by
$$\Delta^n f \colon Q \times V \rightarrow {\R} \colon (q;v) \mapsto \delta^n f(q;v, v, \dots, v).$$
It is the restriction of the $n$-th polarization to the product of the space $Q$ with the diagonal of $V^n$.

\proclaim Definition.
A function $f : Q \to \R$ is said to be a {\it homogeneous polynomial of degree} $n$ at $q\in Q$ if there is a function $F: Q\times V^n \to \R$, multilinear at $q$ in its vectorial arguments, such that, for each $v \in V$,
$$f(q + v) =  {1\over n!} F(q; \underbrace{v,\dots, v}_{n\, {\rm times}}).\eqno\bull$$
\noindent
{\it Polynomials} are sums of homogeneous polynomials.  
The {\it degree} of a polynomial is the highest degree of its non zero homogeneous components.

We denote by $K^n(Q,q)$ the vector space of homogeneous polynomial functions of degree $n$ at $q$.   

\proclaim Proposition.
If $f \in K^n(Q,q)$, then the polarization $\delta^n f$ is multilinear at $q$ and for each $v \in V$,
$$\Delta^n f(q;v)= n! f(q+v). \eqno\bull$$

\proclaim Definition.
The limit
$$ d^n f(q;v_1,v_2,\ldots,v_n) = \lim_{s\to 0}{1\over s^n}\delta^n f (q;sv_1,sv_2, \dots,sv_n), $$
if it exists, is called the $n$-{\it th multidirectional derivative} of $f$ at the point $q\in Q$ in the
 {\it multidirection} $(v_1,v_2,\ldots,v_n) \in V^n$.
 \hfill$\bull$

The $n$-{\it th directional derivative} at $q$ in the direction $v$ is the restriction
$$D^n f(q;v) = d^n f(q;v,v,\ldots,v)$$
of the multidirectional derivative to $Q$ times the diagonal of $V^n$.

\proclaim Proposition.
If $f\in C^{\infty}(Q)$ then for each $i \in N$ and  each $q_0\in Q$, 
$$D^if(q_0, \cdot - q_0) \in K^i(Q,q_0). \eqno\bull$$

\proclaim Theorem (Modified Taylor's theorem). 
If $f\in C^n(Q)$, then for each $q,q_0\in Q$
$$f(q) = \sum_{k=0}^n{1\over k!}D^k f(q_0;q-q_0) + r(q,q_0) \eqno(23)$$
where $r$ satisfies
$$r(q,q)= 0 \quad {\sl and} \quad 
 \lim_{(q,q_0)\to(\bar q,\bar q)} {r(q,q_0) \over \|q-q_0\|^n} = 0$$
 for each $\bar q$.
 \hfill$\bull$

\proclaim Note. 

In the formula (23)  the function $\sum_{k=0}^n{1\over k!}D^k f(q_0;\cdot-q_0)$ is the Taylor polynomial of $f$ of order $n$ at $q'$  and the function $r(\cdot,q')$ is the remainder of $f$ of order $n$ at $q_0$. The latter can be given the usual form of a Lagrange remainder, i.e., the product of a homogeneous polynomial of order $n+1$ at $q_0$ and a suitable function belonging to $C^n(Q)$.

\References
 
 [1] M. Barile, F. Barone and W.M. Tulczyjew, {\it Polarizations and Differential Calculus in Affine Spaces}, to appear in Linear Multilinear Algebra.
 
 [2] D.J. Eck, {\it  Product--Preserving Functors on Smooth Manifolds}, J. Pure Appl. Algebra, {\bf 42}, 133--140 (1986).

[3] C. Ehresmann, {\it  Les prolongements d'une vari\'et\'e diff\'erentiable}, C.R. Acad. Sci., {\bf 233}, 598--600, 777--779 (1951).

[4] J. Gancarzewicz, W. Mikulski and Z. Pogoda, {\it Lifts of Some Tensor Fields and Connections to Product Preserving Functors}, Nagoya Math. J. {\bf 135}, 1--41 (1994).

[5] G. Kainz and P.W. Michor, {\it Natural Transformations in Differential Geometry}, Czech. Math. J., {\bf 37} (112), 584--607 (1987).

 [6] I. Kol\'a\v{r}, {\it On Natural Operators on Vector Fields}, Ann. Global Anal. Geom., {\bf 6}, 109--117 (1988).

 [7] I. Kol\'a\v{r}, P.W. Michor and J. Slov\'ak, {\it Natural Operations in Differential Geometry}, Springer (1993).

 [8] O.O. Luciano, {\it Categories of Multiplicative Functors and Weil's Infinitely Near Points}, Nagoya Math. J., {\bf 109}, 69--89 (1988).
 
 [9] H. Matsumura, {\it Commutative Ring Theory}, Cambridge University Press (1989).
 
 [10] B. Mitchell, {\it  Theory of Categories}, Academic Press, New York and London (1965).
 
 [11] A. Morimoto, {\it  Prolongations of Connections to Tangential Fibre Bundles of Higher Order}, Nagoya Math. J. {\bf 40}, 85--97 (1970).

[12] A. Morimoto, {\it  Prolongations of G-Structure to Tangent Bundles of Higher Order}, Nagoya Math. J. {\bf 38}, 153--179 (1970).

[13] J.  Mu\~noz, F.J. Muriel and J. Rodriguez, {\it  The Contact System on the $(m,\ell)$--Jet Spaces}, Arch. Math. Brno {\bf 37},  291--300 (2001).

[14] J.  Mu\~noz, F.J. Muriel and J. Rodriguez, {\it  On the Finiteness of Differential Invariants}, J. Math. Anal. Appl., {\bf 284}, 266--282 (2003).

[15] A. Weil, {\it Th\'eorie des points proches sur les vari\'et\'es diff\'erentiables}, Colloques internat. Centre nat. Rech. Sci. {\bf 52}, 111--117 (1953).

\bye